\documentclass{elsarticle}
\usepackage[colorlinks,
            linkcolor=blue,
            citecolor=blue,
            urlcolor=magenta,
            linktocpage,
            plainpages=false]{hyperref}
\usepackage{graphicx}
\usepackage{amsmath}
\usepackage{amssymb}
\usepackage{amsfonts}
\usepackage{amsthm}
\usepackage{caption}
\usepackage{subcaption}

\newtheoremstyle{mystyle}
  {}
  {}
  {\itshape}
  {}
  {\bfseries}
  {.}
  { }
  {\thmname{#1}\thmnumber{ #2}\thmnote{ (#3)}}
\theoremstyle{mystyle}

\newtheorem{lemma}{Lemma}
\newtheorem{definition}{Definition}
\newtheorem{proposition}{Proposition}
\newtheorem{corollary}{Corollary}
\newtheorem{theorem}{Theorem}
\newtheorem*{proof*}{Proof}

\newcommand{\romanenumi}{
  \renewcommand{\theenumi}{(\roman{enumi})}
  \renewcommand{\labelenumi}{\theenumi}
}

\newcommand{\alenumi}[1]{
  \renewcommand{\theenumi}{(#1\arabic{enumi})}
  \renewcommand{\labelenumi}{\theenumi}
}

\romanenumi

\newcommand{\klsep}{\operatorname{||}}
\newcommand{\nlim}{\lim_{n \to \infty}}
\newcommand{\reals}{\mathbb{R}}

\newcommand{\mref}[1]{(\ref{#1})}
\newcommand{\dd}{\,d}

\newcommand{\argmax}{\operatornamewithlimits{argmax}}
\newcommand{\indep}{\perp\hspace*{-0.21cm}\perp}
\newcommand{\given}{\operatorname{|}}
\newcommand{\at}{\,\middle|\,}
\newcommand{\Prob}{\operatorname{P}}
\newcommand{\E}{\operatorname{E}}
\newcommand{\VAR}{\operatorname{VAR}}
\newcommand{\A}{\mathbf{A}}
\newcommand{\B}{\mathbf{B}}
\newcommand{\X}{\mathbf{X}}
\newcommand{\Z}{\mathbf{Z}}
\newcommand{\Zi}{\mathbf{Z}_{(i)}}
\newcommand{\ZA}{\mathbf{Z}_{(\A)}}
\newcommand{\ZAB}{\mathbf{Z}_{(\A\B)}}
\newcommand{\G}{\mathcal{G}}
\newcommand{\D}{\mathcal{D}}
\newcommand{\Di}{\D_{(i)}}
\newcommand{\DiG}[1]{\D_{(i)}^{(#1)}}
\newcommand{\PXi}[1]{\mathit{Pa}_{#1}(X_i)}
\newcommand{\PXj}[1]{\mathit{Pa}_{#1}(X_j)}
\newcommand{\XPi}[1]{X_i \given \PXi{#1}}
\newcommand{\T}{\Theta_{X_i}}
\newcommand{\Tz}{\Theta_0}
\newcommand{\hT}{\widehat{\Theta}_{X_i}}
\newcommand{\TT}{\T \given \PXi{}}
\newcommand{\piijk}{\pi_{ik \given j}}
\newcommand{\hpiijk}{\widehat{\pi}_{ik \given j}}
\newcommand{\DXi}{\Delta_{X_i}}
\newcommand{\GXi}{\Gamma_{X_i}}
\newcommand{\DD}{\boldsymbol{\Delta}}
\newcommand{\GG}{\boldsymbol{\Gamma}}
\newcommand{\logl}{\ell}
\newcommand{\nal}{\bar\ell}
\renewcommand{\dim}[1]{\left|#1\right|}
\newcommand{\SGD}{S(\G \given\D)}
\newcommand{\SBIC}{S_{\mathrm{BIC}}(\G \given \D)}
\newcommand{\SML}{S_{\mathrm{ML}}(\G \given \D)}
\newcommand{\SPL}[1]{S_{\mathrm{PL}}(\G_{#1} \given \D)}
\newcommand{\graphs}{\mathbb{G}}
\newcommand{\inG}{\in \graphs}
\newcommand{\KL}{\operatorname{KL}}
\newcommand{\hG}{\widehat{\mathcal{G}}}
\newcommand{\Gz}{\mathcal{G}_0}
\newcommand{\mXd}{\mu_{ij}}
\newcommand{\bXd}{\boldsymbol{\beta}_{ij}}
\newcommand{\hbXd}{\widehat{\boldsymbol{\beta}}_{ij}}
\newcommand{\sXd}{\sigma^2_{ij}}
\newcommand{\hsXd}{\widehat{\sigma}^2_{ij}}
\newcommand{\Ij}{\mathcal{I}_j}
\newcommand{\W}{\mathbf{W}}
\newcommand{\Jac}{\nabla_{\T}}
\newcommand{\Hess}{\mathbf{H}_{\T}\!\!}

\def\bbbone{{\mathchoice {\rm 1\mskip-4mu l} {\rm 1\mskip-4mu l}
            {\rm 1\mskip-4.5mu l} {\rm 1\mskip-5mu l}}}

\begin{document}

\begin{frontmatter}

\title{Learning Bayesian Networks from Incomplete Data with the Node-Average
  Likelihood}

\author[1]{Tjebbe Bodewes}
\address[1]{Zivver, Rotterdam, The Netherlands; formerly
  Department of Statistics, University of Oxford, UK}
\ead{tjebbe.bodewes@linacre.ox.ac.uk}

\author[2]{Marco Scutari}
\address[2]{Istituto Dalle Molle di Studi sull'Intelligenza Artificiale (IDSIA),
     Lugano, Switzerland}
\ead{scutari@idsia.ch}

\begin{abstract}
  Bayesian network (BN) structure learning from complete data has been
  extensively studied in the literature. However, fewer theoretical results are
  available for incomplete data, and most are related to the
  Expectation-Maximisation (EM) algorithm. Balov \cite{balov} proposed an
  alternative approach called Node-Average Likelihood (NAL) that is competitive
  with EM but computationally more efficient; and he proved its consistency and
  model identifiability for discrete BNs.

  In this paper, we give general sufficient conditions for the consistency of
  NAL; and we prove consistency and identifiability for conditional Gaussian
  BNs, which include discrete and Gaussian BNs as special cases. Furthermore, we
  confirm our results and the results in Balov \cite{balov} with an independent
  simulation study. Hence we show that NAL has a much wider applicability than
  originally implied in Balov \cite{balov}, and that it is competitive with EM
  for conditional Gaussian BNs as well.
\end{abstract}

\end{frontmatter}

\section{Background and Notation}
\label{sec:intro}

Bayesian Networks \cite[BNs;][]{pearl,lauritzen} are a class of graphical models
in which the nodes of a directed acyclic graph (DAG) $\G$ represent a set $\X =
\{X_1, \ldots, X_N\}$ of random variables describing some quantities of
interest. The arcs connecting these nodes express direct dependence
relationships, with graphical separation in $\G$ (called \emph{d-separation})
implying conditional independence in probability. Formally, a BN is then defined
as an \emph{independence map} that maps graphical separation to conditional
independence. As a result, $\G$ induces the factorisation
\begin{equation}
  \Prob_{\G}(\X) = \prod_{i=1}^N \Prob\left(\XPi{}, \T\right),
\label{eq:parents}
\end{equation}
in which the joint distribution of $\X$ decomposes into one \emph{local
distribution} for each $X_i$ (with parameters $\T$, $\bigcup_{X_i \in \X} \T =
\Theta$) conditional on its parents $\PXi{}$. Thus BNs provide a compact
and modular representation of both low- and high-dimensional problems.

Multiple DAGs can represent the same set of independence relationships, and can
be grouped into \emph{equivalence classes} \cite{chickering} whose elements are
probabilistically indistinguishable without additional information. All DAGs in
the same equivalence class share the same underlying undirected graph and
v-structures (patterns of arcs like $X_i \rightarrow X_j \leftarrow X_k$, with
no arcs between $X_i$ and $X_k$). It is easy to see that the two other possible
patterns of three nodes and two arcs result in equivalent factorisations:
\begin{multline*}
  \underbrace{
    \Prob(X_i)\Prob(X_j \given X_i)\Prob(X_k \given X_j)
  }_{X_i \rightarrow X_j \rightarrow X_k} =
  \Prob(X_j, X_i)\Prob(X_k \given X_j) = \\
  \underbrace{
    \Prob(X_i \given X_j)\Prob(X_j)\Prob(X_k \given X_j)
  }_{X_i \leftarrow X_j \rightarrow X_k}.
\end{multline*}
Hence each equivalence class is represented by the completed partially-directed
acyclic graph (CPDAG) where only arcs belonging to v-structures and those that
would introduce additional v-structures or cycles are directed.

As for the probability distribution of $\X$, the literature has mostly focused
on three sets of assumptions for analytical and computational reasons.
\emph{Discrete BNs} \cite{heckerman} assume that $\X$ and the $X_i$ are
multinomial random variables. Local distributions take the form
\begin{align}
  &\XPi{} \sim \mathit{Mul}\left(\piijk\right),&
  &\piijk = \Prob\left(X_i = k \given \PXi{} = j\right);
\label{eq:cpt}
\end{align}
their parameters $\piijk$ are the conditional probabilities of $X_i$ given each
configuration of the values of its parents, usually represented as a conditional
probability table for each $X_i$.

\emph{Gaussian BNs} \citep[GBNs;][]{heckerman3} model $\X$ with a multivariate
normal random variable and assume that the $X_i$ are univariate normals linked
by linear dependencies. The parameters of the local distributions can be
equivalently written \citep{weatherburn} as the partial correlations
$\rho_{X_i, X_j \given \PXi{} \setminus X_j}$ between $X_i$ and each parent
$X_j$ given the other parents; or as the coefficients $\boldsymbol{\beta}_{X_i}$
of the linear regression model
\begin{align}
  &X_i = \mu_{X_i} + \PXi{}\boldsymbol{\beta}_{X_i} + \varepsilon_{X_i},&
  &\varepsilon_{X_i} \sim N\left(0, \sigma^2_{X_i}\right), \sigma^2_{X_i} > 0,
\label{eq:reg}
\end{align}
so that $\XPi{} \sim N\left(\mu_{X_i} + \PXi{}\boldsymbol{\beta}_{X_i},
\sigma^2_{X_i}\right)$.

Finally, \emph{conditional linear Gaussian BNs} \citep[CLGBNs;][]{unification}
combine discrete and continuous random variables in a mixture-of-Gaussians model
as follows:
\begin{itemize}
  \item Discrete $X_i$ are only allowed to have discrete parents (denoted
    $\DXi$), and are assumed to follow a multinomial distribution as in
    \mref{eq:cpt}.
  \item Continuous $X_i$ are allowed to have both discrete and continuous
    parents (denoted $\GXi$, $\DXi \cup \GXi = \PXi{}$), and their local
    distributions are
    \begin{equation*}
      \XPi{} \sim N\left(\mXd + \GXi\bXd, \sXd\right)
    \end{equation*}
    or, equivalently,
    \begin{align}
      &X_i = \mXd + \GXi\bXd + \varepsilon,&
      &\varepsilon \sim N(0, \sXd), \sXd > 0,
    \label{eq:cgnode}
    \end{align}
    which define a mixture of linear regressions against the continuous parents
    $\GXi$ with one component for each configuration $j$ of the discrete parents
    $\DXi$. Hence $\T = \{ \mXd, \bXd, \sXd\}$. If $X_i$ has no discrete
    parents ($\DXi = \varnothing$), the mixture in \mref{eq:cgnode} simplifies
    to a single linear regression like that in \mref{eq:reg}.
\end{itemize}
We denote discrete nodes with $\DD$ and Gaussian nodes as $\GG$, with $\DD \cup
\GG = \X$. If $\DD = \varnothing$ or $\GG = \varnothing$, a CGBN reduces
respectively to a GBN or to a discrete BN. For this reason, we will derive our
results only for CGBNs but they will hold for both discrete BNs and GBNs as
particular cases.

Other distributional assumptions have seen less widespread adoption for various
reasons: for instance, copulas \cite{copula} and truncated exponentials
\cite{truncexp} lack exact conditional inference and simple closed-form
estimators.

The task of learning a BN from a data set $\D$ is performed in two steps:
\begin{equation}
  \underbrace{\Prob(\G, \Theta \given \D)}_{\text{learning}} =
    \underbrace{\Prob(\G \given \D)}_{\text{structure learning}} \cdot
    \underbrace{\Prob(\Theta \given \G, \D)}_{\text{parameter learning}}.
\label{eq:lproc}
\end{equation}
\emph{Structure learning} consists in finding the DAG $\G$ that
encodes the dependence structure of the data in the space $\graphs$ of all
possible DAGs, thus maximising $\Prob(\G \given \D)$ or some alternative
goodness-of-fit measure. \emph{Parameter learning} consists in estimating the
parameters $\Theta$ given the $\G$ obtained from structure learning, that is
$\argmax_\Theta \Prob(\Theta \given \G, \D)$. If we assume that parameters in
different local distributions are independent and that data contain no missing
values, the $\T$ can be learned separately for each node because
\mref{eq:parents} implies
\begin{equation*}
  \Prob(\Theta \given \G, \D) = \prod_{i=1}^N \Prob(\TT, \D).
\end{equation*}
Assuming that $\G$ is sparse, each $X_i$ will have a small number of parents and
therefore $\Prob(\TT, \D)$ will be a low-dimensional distribution that can be
estimated efficiently from the data.

On the other hand, structure learning is NP-hard \cite{nphard}. Many algorithms
have been proposed for this problem following one of three possible approaches:
\emph{constraint-based}, \emph{score-based} and \emph{hybrid}. Constraint-based
algorithms are based on the seminal work of Pearl on causal graphical models
\cite{ic}: they use conditional independence tests to assess which conditional
independence relationships are supported by data, and they assume that $\G$ is a
perfect map of to identify which arcs to include in $\G$ from the results of
those tests. The most popular algorithm in this class is arguably the PC
algorithm \cite{colombo}, which was originally introduced in \cite{spirtes} in
the context of causal discovery. Score-based algorithms represent the
application of general-purpose optimisation techniques to BN structure learning:
each candidate DAG is assigned a network score reflecting its goodness of fit,
which the algorithm then attempts to maximise. Some examples are heuristics such
as greedy search, simulated annealing \cite{bouckaert} and genetic algorithms
\cite{larranaga}; they work well in practical applications but do not come with
any formal guarantees about the optimality of the networks they learn. Other
algorithms \cite{order1,k2}, such as K2, seek to speed-up structure learning by
first obtaining a topological ordering of the nodes and then learning the
optimal DAG conditional on that. Exact maximisation algorithms have also become
feasible in recent years thanks to increasingly efficient pruning of the space
of DAGs and tight bounds on the scores \cite{cutting,suzuki17}. Finally, hybrid
algorithms combine the previous two approaches, using conditional independence
tests to reduce the number of candidate DAGs to score while searching for the
DAG maximising some network score. The best-known member of this family is the
Max-Min Hill Climbing algorithm by \cite{mmhc}. A comprehensive review and
comparison of these three approaches for complete data can be found in
\cite{ijar19}; in the following, we will only consider score-based structure
learning.

Both parameter and structure learning are mathematically and computationally
harder problems when learning BNs from incomplete data because decomposition in
\mref{eq:parents} no longer holds. A common approach is to use the
Expectation-Maximisation algorithm \cite[EM;][]{em} for both parameter learning
and structure learning to reuse methods originally proposed for complete data.
However, this choice comes at a significant computational cost. This is
particularly the case for structure learning, where complete-data algorithms are
embedded within EM in the Structural EM \cite{friedman97,friedman98} algorithm.
To address this issue, Balov \cite{balov} proposed an alternative approach based
on the Node-Average Likelihood (NAL) that is competitive with EM-based
approaches in terms of structural accuracy but has a much lower computational
cost. He then proved both identifiability and the consistency of NAL for
discrete BNs. Both approaches, as well as others, are reviewed in
\cite{neerlandica19} which provides a broad overview of parameter and structure
learning from incomplete data.

In this paper we establish general conditions for both properties of NAL
and we show that they hold for CGBNs. In Section \ref{sec:complete} we
briefly review score-based learning from complete data, moving to incomplete
data in Section \ref{sec:incomplete}. Our novel results and the required
regularity conditions are introduced in Section \ref{sec:theorems} for both
identifiability (Section \ref{sec:identifiability}) and consistency (Section
\ref{sec:consistency}). Proofs for all but the main consistency result are
included in \ref{sec:proofs}. Finally, in Section \ref{sec:experiments} we
confirm the results on discrete BNs obtained by Balov \cite{balov} with a
simulation study, and that similar results hold for CGBNs.

\section{Structure Learning from Complete Data}
\label{sec:complete}

Score-based BN structure learning consists of two components: a score function
$\SGD$ and an algorithm that determines how we explore the space of the DAGs.
Each candidate DAG is assigned a score $\SGD$ reflecting its goodness of fit,
which the algorithm then attempts to maximise to obtain a (possibly local)
optimum DAG as $\argmax_{\G \inG} \SGD$. Heuristic algorithms such as tabu
search \cite{glover} are more common in practical applications, but exact
algorithms have been proposed as well \cite{cutting,suzuki17} although their use
is limited to small and medium problems by their computational complexity. As
for the score function, using Bayes theorem we have that
\begin{equation*}
  \Prob(\G \given \D) \propto \Prob(\G)\Prob(\D \given \G) =
    \Prob(\G)\int \Prob(\D \given \G, \Theta) \Prob(\Theta \given \G) \dd\Theta
\label{eq:marginal-likelihood}
\end{equation*}
where $\Prob(\G)$ is the prior distribution over the space of the DAGs spanning
the variables in $\X$ and $\Prob(\D \given \G)$ is the marginal likelihood of
the data given $\G$ averaged over all possible parameter sets $\Theta$.
The most common choice in the literature for $\Prob(\G)$ is the uniform
distribution. The space of the DAGs grows super-exponentially in $N$
\citep{harary} and that makes it cumbersome to specify granular priors that
assign different prior probabilities to individual DAGs. However, informative
priors that assign probabilities to groups of DAGs based on their structural
properties have been proposed in the literature, including the sparsity-inducing
priors presented in \cite{pgm16} and \cite{gprior} and the informative priors
from \cite{csprior} and \cite{mukherjee}.

From \mref{eq:parents}, the marginal likelihood decomposes into one term for
each node,
\begin{equation}
  \SML = \prod_{i=1}^N \int \Prob(\XPi{}, \T) \Prob(\TT) \dd\T,
\label{eq:marginal}
\end{equation}
which makes it possible to use it as a network score (denoted $\SML$) in a
computationally efficient way: the fact that $\SML$ is decomposable speeds up
structure learning considerably because \mref{eq:marginal} is
available in closed form for discrete BNs, GBNs and CGBNs \cite{unification}
and because we only recompute differing portions of $\Prob(\D \given \G)$ as we
score and compare DAGs. It is also consistent for complete data.

Due to the difficulty of choosing a prior over $\Theta$ and the resulting
performance implications \cite{ambn17}, a common alternative is the
\emph{Bayesian Information Criterion} \cite[BIC;][]{schwarz},
\begin{align}
  \logl(\G, \Theta \given \D)
    &= \frac{1}{n}\sum\nolimits_{X_i \in \X} \log\Prob(\XPi{}, \hT),
  \label{eq:loglik} \\
  \SBIC &= \logl(\G, \Theta \given \D) - \frac{\log(n)}{2n}\dim{\Theta};
  \label{eq:sbic}
\end{align}
where $n = |\D|$ is the sample size and the $\hT$ are the maximum likelihood
estimates (MLEs) of the $\T$ for $\D$. In contrast, we will denote by $\logl(\G,
\Theta)$ and $\logl(\XPi{}, \T)$ the population log-likelihoods. $\SBIC$ is a
particular case of the penalised log-likelihood
\begin{align}
  &S_{\mathrm{PL}}(\G \given \D) = \logl(\G, \Theta \given \D) - \lambda_n h(\G),&
  &\lambda_n \geq 0
\end{align}
where $\lambda_n$ is a penalisation coefficient and $h: \graphs \to \reals^+$
measures model complexity such that $\G_1 \subset \G_2 \Rightarrow
h(\G_1) < h(\G_2)$.\footnote{We say $\G_1 \subseteq \G_2$ if the arc set of
$\G_1$ is a subset of that of $\G_2$.} $S_{\mathrm{PL}}(\G \given \D)$ is
related to the $\ell_0$ penalised maximum likelihood studied in \cite{l0} for
GBN structure learning, which has the same general form; but it does not make
any distributional assumption. $\SBIC$ is decomposable, since $h(\G) =
\dim{\Theta} = \sum_{\X} \dim{\T}$ which is the number of parameters of each
local distribution; and it does not depend on any hyperparameter. Furthermore,
it is equivalent to the minimum description length \cite{rissanen} of $\G$ and
it is asymptotically equivalent to $\SML$. Hence, it is consistent for complete
data. Setting $\lambda_n = 1/n$ instead of $\lambda_n = \log(n)/2n$ gives the
\emph{Akaike Information Criterion} \cite[AIC;][]{akaike} which, on the other
hand, is not consistent for complete data \cite{bozdogan}.

An important property shared by AIC and BIC is \textit{score equivalence}
\cite{chickering}: they take the same value for all DAGs in the same equivalence
class. The same is true for $\SML$ for particular choices of the prior over
$\Theta$, as shown in \cite{unification}. This acknowledges that those DAGs are
probabilistically indistinguishable from the data, and thus cannot be ranked by
the score function. We can, however, still express a preference for some DAGs
over others in the same equivalence class by means of informative priors
$\Prob(\G)$ over the space of DAGs.

\section{Structure Learning from Incomplete Data}
\label{sec:incomplete}

In the context of BNs, incomplete data are modelled using auxiliary nodes $Z_i
\in \Z$ that encode whether the corresponding $X_i$ is observed for each
observation. (Say, $Z_i = 1$ if $X_i$ is observed and $Z_i = 0$ if not.) Rubin
and Little \cite{rubin2,rubin1} formalised three possible patterns (or
mechanisms) of missingness:
\begin{itemize}
  \item \emph{Missing completely at random} (MCAR): when complete samples are
    indistinguishable from incomplete ones. In other words, the probability that
    a value will be missing is independent from both observed and missing
    values.
  \item \emph{Missing at random} (MAR): cases with incomplete data differ from
    cases with complete data, but the pattern of missingness is predictable from
    other observed variables. In other words, the probability that a value will
    be missing is a function of the observed values.
  \item \emph{Missing not at random} (MNAR): the pattern of missingness is not
    random or it is not predictable from other observed variables; the
    probability that an entry will be missing depends on both observed and
    missing values. Common examples are variables that are missing
    systematically or for which the patterns of missingness depend on the
    missing values themselves.
\end{itemize}
These patterns of missingness can be modelled graphically: MCAR implies $\Z
\indep \X$; MAR implies $\Z \indep \X$ for the incomplete observations
conditional on the complete observations; and MNAR does not imply any
independence constraint. Note that, however, in the following we will still
consider $\G$ to be spanning just $\X$ because we will restrict ourselves to the
MCAR case in which the dependencies between the $\X$ and the $\Z$ are completely
determined.

When the data are incomplete, scoring $\G$ requires the missing values to be
integrated out. We can make this apparent by rewriting $\Prob(\D \given \G)$ as
a function of $\D^O$ (the observed data) and $\D^M$ (the missing data):
\begin{equation*}
  \Prob\left(\D \given \G\right) = \Prob\left(\D^O, \D^M \given \G\right)
    = \int \Prob\left(\D^O, \D^M \given \G, \Theta\right)
           \Prob\left( \Theta \given \G \right) \dd\Theta.
\end{equation*}
Averaging over all the possible configurations of $\D^M$ as well leads to
\begin{equation}
  \Prob\left(\D^O \given \G\right)
    = \iint \Prob\left(\D^M \given \D^O, \G, \Theta\right)
            \Prob\left(\D^O \given \G, \Theta\right)
            \Prob\left(\Theta \given \G\right) \dd\Theta \dd\D^M
\label{eq:marginal2}
\end{equation}
which contains one extra dimension for each missing value (in addition to one
dimension for each parameter in $\Theta$), making it infeasible to compute in
practical applications. An additional problem is that, while $\Prob(\D^O \given
\G, \Theta)$ decomposes as in \mref{eq:marginal}, $\Prob(\D^M \given \D^O, \G,
\Theta)$ does not in the general case.

An alternative is to use the Expectation-Maximisation algorithm
\cite[EM;][]{em}, an iterative procedure consisting of an E(xpectation)-step and
a M(aximisation)-step. In the E-step, we compute the expected sufficient
statistics conditional on the observed data using belief propagation
\cite{pearl,grem,shafer}. In the M-step, complete-data learning methods can be
applied using the expected sufficient statistics instead of the (unobservable)
empirical ones.

There are two ways to apply EM to structure learning. Firstly, we can apply EM
separately to each candidate DAG to be scored, as in the variational Bayes EM
\cite{variational-em}. However, structure learning often involves many
evaluations of the score function, thus making this approach computationally
infeasible beyond small-scale problems. Secondly, we can embed structure
learning in the M-step as follows:
\begin{itemize}
  \item in the E-step, we complete the data by computing the expected sufficient
    statistics using the current network structure;
  \item in the M-step, we find the structure that maximises the expected score
    function for the completed data.
\end{itemize}
This approach is called Structural EM, and has been implemented using both
$\SBIC$ \cite{friedman97} and $\SML$ \cite{friedman98}. While faster than
variational Bayes EM, Structural EM is still computationally expensive.
\cite{friedman97} notes: ``Most of the running time during the execution our
procedure is spent in the computations of expected statistics. This is where our
procedure differs from parametric EM. In parametric EM, we know in advance which
expected statistics are required. [\ldots] In our procedure, we cannot determine
in advance which statistics will be required. Thus, we have to handle each query
separately. Moreover, when we consider different structures, we are bound to
evaluate a larger number of queries than parametric EM.''

Balov \cite{balov} proposed a more scalable approach for discrete BNs under MCAR
called Node-Average Likelihood (NAL). While Balov \cite{balov} defined NAL
relying on the specific form of the multinomial log-likelihood, we will present
it here using a more general notation that allows its extension to CGBNs.
Starting from \mref{eq:parents}, he proposed to compute each term using the $\Di
\subseteq \D$ locally-complete data for which $X_i, \PXi{}$ are observed (that
is, $\Zi = \Z_{X_i, \PXi{}} = 1$):
\begin{equation*}
  \nal(\XPi{}, \hT) =
    \frac{1}{|\Di|} \sum\nolimits_{\Di}\log \Prob(\XPi{}, \hT),
\end{equation*}
which is an empirical estimate of the average node log-likelihood
$\E[\logl(\XPi{})]$. Replacing \mref{eq:loglik} with the above gives
\begin{equation*}
  \nal(\G, \Theta \given \D) = \sum\nolimits_{X_i \in \X} \nal(\XPi{}, \hT)
\end{equation*}
which Balov \cite{balov} used to redefine the penalised log-likelihood score
function as
\begin{equation}
  \SPL{} = \nal(\G, \Theta \given \D) - \lambda_n h(\G)
  \label{eq:spl}
\end{equation}
and structure learning as $\hG = \argmax_{\G \inG} \SPL{}$. NAL makes a more
efficient use of incomplete data than discarding all incomplete samples; and it
does not incur the computational costs of EM approaches since it can be
plugged in complete-data structure learning approaches. However, comparing two
DAGs means that NAL is evaluated on potentially different subsets of $\D$ for
each $X_i$ in different DAGs; hence the usual results on MLEs and nested models
do not apply. Balov \cite{balov} proved both identifiability and consistency of
score-based structure learning when using $\SPL{}$ for discrete BNs. We will now
prove both properties hold more generally, and in particular that they hold for
CGBNs.

\section{Properties of Node-Average Likelihood}
\label{sec:theorems}

We study two properties of NAL: whether the true DAG $\Gz$ is
\emph{identifiable}, and under which conditions $\hG$ is a \emph{consistent}
estimator of $\Gz$. For each of our results, we reference the corresponding
theorems that Balov \cite{balov} derived for discrete BNs. Our contribution is
the generalisation of these results to CGBNs under mild regularity conditions,
which requires completely different proofs since Balov \cite{balov} relied
materially on the form of the multinomial log-likelihood and on not having both
continuous and discrete parents in each $\PXi{}$.

\subsection{Identifiability}
\label{sec:identifiability}

In this section we show under what conditions NAL can be used to identify the
true $\Gz$ of the BN from data in the limit of the sample size. Firstly, we
prove that NAL is non-decreasing in the size of the parent sets and thus
overfits like the log-likelihood does for complete data. This motivates the use
of $\SPL{}$ to avoid favouring overly complex models.
\begin{lemma}[L7.1 in \cite{balov}]
  For any $X_i$ and disjoint $\A, \B \subset \X$ such that \linebreak $\Z_\B \indep \{X_i,
  \A, \B\} \given \Z_{X_i, \A} = 1$, we have $\bar\logl(X_i \given \A) \leqslant
  \bar\logl(X_i \given \A, \B)$, with equality if and only if
  $X_i \indep \B \given (\A, \Z_{X_i, \A} = 1)$.
  \label{lemma:nondec}
\end{lemma}
Lemma \ref{lemma:nondec} allows us to state that if $\Gz$ is identifiable, we
can learn it by finding the simplest DAG that maximises NAL. We characterise
this fact in a way that is more appropriate for BN structure learning in the
same way as Balov \cite{balov} below.

\begin{definition}[D3.1 in \cite{balov}]
  $\Gz$ is identifiable if for any $\G \inG$ we have
  $\nal(\G, \Theta) \leqslant \nal(\Gz, \Tz)$ when $\Gz \subseteq \G$;
  and $\nal(\G, \Theta) < \nal(\Gz, \Tz)$ when $\Gz \not\subseteq \G$.
  \label{defIdentifiability}
\end{definition}

We next establish under which conditions $\Gz$ is guaranteed to be identifiable:
under MCAR $\nal(\G, \Theta)$ attains its maximum at all $\G \supseteq \Gz$, and
all these DAGs induce the true distribution
$\Prob_{\Gz}(\X) = \prod_\X \Prob(\XPi{\Gz}, \Theta_{0, X_i})$.
\begin{proposition}[P3.1 in \cite{balov}]
  Under MCAR, we have:
  \begin{enumerate}
    \item
      $\max_{\G \inG} \nal(\G, \Theta) = \nal(\Gz, \Tz)$.
      \label{part:one}
    \item if $\nal(\G, \Theta) = \nal(\Gz, \Tz)$, then
      $\Prob_{\G}(\X) = \Prob_{\Gz}(\X)$. \label{part:two}
    \item if $\Gz \subseteq \G$, then
      $\nal(\G, \Theta) = \nal(\Gz, \Tz)$. \label{part:three}
  \end{enumerate}
  \label{prop:id}
\end{proposition}
The identifiability of $\Gz$ up to its equivalence class $[\Gz]$ follows from
the above and is formally stated below.
\begin{corollary}[C3.2]
  $[\Gz]$ is identifiable under MCAR, that is
  \begin{equation*}
    \Gz \cong \min \left\{\G_* \inG: \bar\ell(\G_*, \Theta_*) =
      \max_{\G \inG} \nal(\G, \Theta) \right\}.
  \end{equation*}
  \label{cor:equiv}
\end{corollary}

\subsection{Consistency}
\label{sec:consistency}

In this section we show that the candidate $\hG$ chosen by maximising $\SPL{}$
is a consistent estimator of the true $\Gz$ under MCAR, and the mild regularity
conditions this requires.

\subsubsection{Regularity Conditions}
\label{sec:regularity}

For all $\G \inG$ and $X_i \in \X$:
\alenumi{R}
\addtolength{\leftmargini}{1.5\labelsep}
\begin{enumerate}
  \item $\hT$ must exist, converge in probability to the population $\T$ with
    speed $O(n^{-1/2})$, and make $\Jac\nal(\XPi{}, \hT)$ vanish.\label{cond:mle}%
  \item The Hessian $\Hess\left(\nal(\XPi{}, \T)\right)$ exists and has finite
    expectation.\label{cond:secder}%
  \item The variance $\nu^2 = \VAR\left(\nal(\XPi{}, \T)\right) < \infty$.\label{cond:var}%
\end{enumerate}
\addtolength{\leftmargini}{-1.5\labelsep}
\romanenumi
Conditions \ref{cond:mle} and \ref{cond:secder} imply that the NAL evaluated at
the MLE is close to the NAL at the population $\T$ with high probability for
large $n$; \ref{cond:var} allows the use of various Central Limit Theorems
\cite[CLTs;][]{billingsley}. \ref{cond:mle}, \ref{cond:secder}, \ref{cond:var}
hold for CGBNs.
\begin{proposition}
  Consider a CGBN over $\X$. If  $\T > 0$ for all $X_i \in \DD$, and $\sXd > 0$
  for all $X_i \in \GG$ and all values $j$ of $\DXi$, then \ref{cond:mle},
  \ref{cond:secder}, \ref{cond:var} are satisfied.
  \label{prop:regCGBN}
\end{proposition}

A more general approach would be to model the $X_i \in \GG$ with a mixture of
generalised linear models \cite{glm}. When using canonical link functions,
\ref{cond:mle} and \ref{cond:secder} simplify due to the results in
\cite{fahrmeir}. However, verifying \ref{cond:mle}, \ref{cond:secder} and
\ref{cond:var} is non-trivial hence we leave this as a potentially tractable
case for future work.

\subsubsection{Consistency Results}
\label{sec:results}

Starting from the consistency of $\nal(\G, \Theta \given \D)$, we will now
establish for which sequences of $\lambda_n$ the $\SPL{}$ in \mref{eq:spl} is
consistent for complete data and under MCAR, that is, $\nlim P(\hG = \Gz) = 1$.
We note the sufficient conditions for the consistency of $\hG$ from
Balov \cite[P3.2, C4.1]{balov}:
\alenumi{C}
\addtolength{\leftmargini}{1.5\labelsep}
\begin{enumerate}
  \item If $\Gz \subseteq \G_1$ and $\Gz \not\subseteq \G_2$, then
    $\nlim \Prob(\SPL{1} > \SPL{2}) = 1$.\label{suff:one}%
  \item If $\Gz \subseteq \G_1$ and $\G_1 \subset \G_2$, then
    $\nlim \Prob(\SPL{1} > \SPL{2})= 1$.\label{suff:two}%
  \item $\exists\, \G, X_i, X_j$, $i \neq j$ such
    that $\PXi{\Gz} \subset \PXi{\G}$, \linebreak $\PXj{\G} = \PXj{\Gz}$
    and $\Prob(\Z_{\PXi{\G} \setminus \PXi{\Gz}} \given
    \Zi^{(\Gz)}) \in (0,1)$.\label{suff:three}%
\end{enumerate}
\addtolength{\leftmargini}{-1.5\labelsep}
\romanenumi
These conditions can be trivially extended to equivalence classes as in
Balov \cite[C3.3]{balov}.

To prove consistency, we first establish some intermediate results. Unlike Balov
\cite{balov}, we need a rigorous treatment of scoring DAGs that may represent
misspecified models \cite{white} that are not representable in terms of $\Gz$;
an example would be a $X_i \in \GG$ being a mixture of regressions in $\G$ and a
single linear regressions in $\Gz$.

Lemma \ref{lemma:toscore} shows that if the difference in NAL between two
DAGs is $O(n^{-\alpha})$, then the less complex DAG is chosen asymptotically if
$\lambda_n \to 0$ slower than $n^{-\alpha}$.
\begin{lemma}
  For $\G_1, \G_2 \inG$, if
  \begin{enumerate}
    \item $\nal(\G_1, \Theta_1 \given \D) -
      \nal(\G_2, \Theta_2 \given \D) = O(n^{-\alpha})$ for $\alpha > 0$;
    \item $h(\G_1) < h(\G_2)$;
    \item $\nlim n^{\alpha}\lambda_n = \infty$;
  \end{enumerate}
  then $\nlim P(\SPL{1} > \SPL{2}) = 1$.
\label{lemma:toscore}
\end{lemma}

Lemma \ref{lemma:mle} establishes that the difference between the NAL at $\hT$
and the NAL at the population $\T$ is $O(n^{-1})$, which is relevant due to
Lemma \ref{lemma:toscore}. Working with the NAL at the true
parameters allows us to exploit conditional independencies and turns the NAL
into a sum of i.i.d. random variables, allowing the Lindeberg-L\'evy and the
Lindeberg-Feller \cite[p.~369]{billingsley} CLTs to be applied.
\begin{lemma}
  If \ref{cond:mle} and \ref{cond:secder} hold, then for all $\G$ and $X_i$:
  \begin{equation*}
    \nal(\XPi{}, \hT) - \nal(\XPi{}, \T) = O(n^{-1}).
  \end{equation*}
\label{lemma:mle}
\end{lemma}

Lemma \ref{lemma:cpd} establishes that if $\G \supseteq \Gz$ the local
distributions for $\G$ reduce to those of $\Gz$, and that population and sample
distributions coincide. These results are crucial in linking $\nal(\G, \Theta
\given \D)$ and $\nal(\G, \Theta)$ in Lemma \ref{lemma:lconv}.
\begin{lemma}
  If $\G \supseteq \Gz$, then for each $X_i \in \X$:.
\begin{enumerate}
    \item $\Prob(\XPi{\G}, \T, \Zi) = \Prob(\XPi{\Gz}, \T, \Zi)$.
    \item $\Prob(\XPi{\G}, \T, \D) = \Prob(\XPi{\Gz}, \T, \D)$ almost
      surely.
    \item $\Prob(\XPi{\G}, \T, \D) = \Prob(\XPi{\G}, \T, \Zi)$ almost
      surely.
  \end{enumerate}
\label{lemma:cpd}
\end{lemma}

Lemma \ref{lemma:lconv} links the sample NAL $\nal(\G, \Theta \given \D)$ and
the population NAL $\nal(\G, \Theta)$ and establishes the convergence rate of
the former. It will be used to show \ref{suff:one} in Theorem \ref{theorem:main}
by way of Proposition \ref{prop:id}.

\begin{lemma}[L4.1 in \cite{balov}]
  If \ref{cond:mle}, \ref{cond:secder}, \ref{cond:var} are satisfied, then
  \begin{align}
    &\nlim \Prob(\nal(\G, \Theta \given \D) \leqslant \nal(\G, \Theta)) = 1&
    &\text{for any $\G \inG$.}
    \label{eq:negativity}
  \end{align}
  Furthermore, $\nal(\G, \Theta \given \D) - \nal(\G, \Theta) =
  O(n^{-1/2})$ if $\Gz \subseteq \G$.
\label{lemma:lconv}
\end{lemma}

Theorem \ref{theorem:main} is the key result of this section, showing that
BIC is consistent for complete data but it is not under MCAR. AIC is not
consistent in either case.

\begin{theorem}[T4.1 in \cite{balov}]
  Let $\Gz$ be identifiable, $\lambda_n \to 0$ as $n \to \infty$, and assume
  that \ref{cond:mle}, \ref{cond:secder} are satisfied. Then as $n \to \infty$:
  \begin{enumerate}
    \item If $\Prob(\Z = 1) = 1$ and $n\lambda_n \to \infty$, $\hG$ is
      consistent.
    \item Under MCAR and \ref{cond:var}, if $\sqrt{n} \lambda_n \to \infty$,
      $\hG$ is consistent.
    \item Under MCAR, \ref{suff:three} and \ref{cond:var}, if
      $\lim_{n\to \infty} \sqrt{n}\lambda_n < \infty$, then $\hG$ is
      not consistent.
  \end{enumerate}
\label{theorem:main}
\end{theorem}
\begin{proof*}[Theorem \ref{theorem:main}]
  For $\G_1, \G_2 \inG$ such that $\Gz \subseteq \G_1$, $\Gz \not \subseteq
  \G_2$, we must show that $ \nlim \Prob(\SPL{1} > \SPL{2}) = 1$.

  Lemma \ref{lemma:nondec} implies that $\nal(\G_1, \Theta_1) \geqslant \nal(\Gz,
  \Tz)$, while identifiability gives \linebreak $\nal(\G_2, \Theta_2) < \nal(\Gz,
  \Tz)$, resulting in $\nal(\G_1, \Theta_1) > \nal(\G_2, \Theta_2)$. Lemma
  \ref{lemma:lconv} implies that $\nal(\G_1, \Theta_1 \given \D) \to
  \nal(\G_1, \Theta_1)$ and $\Prob(\nal(\G_2, \Theta_2 \given \D)
  \leqslant \nal(\G_2, \Theta_2)) = 1$ as $n \to \infty$. We deduce that
  $\nlim \Prob(\nal(\G_1, \Theta_1 \given \D) > \nal(\G_2, \Theta_2 \given \D)) = 1$.
  As $\lambda_n \to 0$ for $n \to \infty$, we have
  $\nlim \Prob(\SPL{1} > \SPL{2}) = 1$
  and condition \ref{suff:one} is satisfied for parts (i)-(iii).

  We show that under the conditions in parts (i) and (ii),
  \begin{equation*}
    \nlim \Prob(\SPL{1}) > \SPL{2}) = 1
  \end{equation*}
  for $\G_1, \G_2 \inG$ such that $\Gz \subseteq \G_1 \subset \G_2$.
  This means that \ref{suff:two} is satisfied and $\hG$
  is consistent. We then show that under the assumptions in part (iii), there
  exists a $\G \supset \Gz$ such that $\nlim \Prob(\SPL{} > \SPL{0}) > 0$, which
  implies inconsistency of $\hG$. Note that if $\Prob(\Z = 1) = 1$ as in part
  (i), then $\X \indep \Z$. Therefore, results derived under MCAR hold for all
  three parts.

  Let $\G \supseteq \Gz$. In the sample NAL, we first apply Lemma
  \ref{lemma:mle} to replace the MLEs with their population values and use
  Lemma \ref{lemma:cpd}(ii) to eliminate the redundant parents:
  \begin{align*}
    \nal(\XPi{\G}, \hT)
      &= \nal(\XPi{\G}, \T) + O(n^{-1}) \\
      &= \frac{1}{|\DiG{\G}|} \sum_{\DiG{\G}} \log \Prob(\XPi{\Gz}, \T)
            + O(n^{-1}).
  \end{align*}
  The difference between the last expression and $\nal(\XPi{\Gz}, \T)$, in
  which $\DiG{\Gz}$ would appear instead of $\DiG{\G}$, represents the
  difference in NAL due to $\Zi^{(\G)} \neq \Zi^{(\Gz)}$. We denote it as
  $d(\Di) + O(n^{-1})$ in the following; the contrast between part (i) and parts
  (ii), (iii) follows from its behaviour for complete and MCAR data.

  \emph{Part (i):}
  For all $\G \supseteq \Gz$, $\nal(\G, \Theta \given \D) - \nal(\Gz, \Tz
  \given \D) = O(n^{-1})$. As $\Gz \subseteq \G_1 \subset \G_2$, it follows that
  $\nal(\G_1, \Theta_1 \given \D) - \nal(\G_2, \Theta_2 \given \D) = O(n^{-1})$
  and $h(\G_1) < h(\G_2)$. From Lemma \ref{lemma:toscore},
  \begin{equation*}
    \nlim \Prob(\SPL{1} > \SPL{2}) = 1.
  \end{equation*}
  Hence \ref{suff:two} is satisfied and part (i) follows.

  \emph{Part (ii):}
  Consider $\G \supset \Gz$: $\PXi{\G} \supset \PXi{\Gz}$ for at least one $X_i$
  and $\Zi^{(\G)}$, $\Zi^{(\Gz)}$ are potentially different, hence
  $\Prob(d(\Di) \neq 0) > 0$. Then by \ref{cond:var} and Balov \cite[L7.4]{balov}:
  \begin{equation*}
    \sqrt{n}d(\Di) \to N(0, \gamma\nu^2),
  \end{equation*}
  where $\gamma = \left(1 - \Prob(\Zi^{(\G)} \given \Zi^{(\Gz)})\right)/
  \Prob(\Zi^{(\G)} \given \Zi^{(\Gz)})$. Reintroducing the $O(n^{-1})$ term and
  noting that $\sqrt{n} O(n^{-1}) \to 0$,
  \begin{equation}
    \sqrt{n}\left(\nal(\XPi{\G}, \hT) - \nal(\XPi{\Gz}, \hT)\right)
      \to N(0, \gamma\nu^2).
  \label{eq:normconv}
  \end{equation}
  As $N(0, \gamma\nu^2)$ is $O(1)$, we obtain that
  $\nal(\G, \Theta \given \D) - \nal(\Gz, \Tz \given \D) = O(n^{-1/2})$
  for all $\G \supset \G_0$.
  It follows that $\nal(\G_1, \Theta_1 \given \D) -
  \nal(\G_2, \Theta_2 \given \D) = O(n^{-1/2})$ for all $\G_1, \G_2 \inG$ such
  that $\Gz \subseteq \G_1 \subset \G_2$.
  Applying Lemma \ref{lemma:toscore} as in part (i), we obtain that if
  $\nlim \sqrt{n}\lambda_n = \infty$, then
  \begin{equation*}
    \limsup_{n\to \infty} \Prob(\SPL{1} > \SPL{2}) = 1.
  \end{equation*}
  Hence \ref{suff:two} is satisfied and part (ii) follows.

  \emph{Part (iii):}
  By Condition \ref{suff:three}, there exist a $\G \inG$ and an $X_i$ such that
  $\PXj{\G} = \PXj{\Gz}$ for all $j \neq i$ and
  $\PXi{\G} \supset \PXi{\Gz}$. Let $h(\G) - h(\Gz) = c > 0$. Then
  $\SPL{} > \SPL{0}$ implies:
  \begin{align}
    \sqrt{n} &\left(\nal(\G, \Theta \given \D) -
         \nal(\Gz, \Tz \given \D)\right) >
              \sqrt{n}\lambda_n \left(h(\G) - h(\Gz)\right), \label{eq:toprove} \\
    \sqrt{n} &\left(\nal(\XPi{\G}, \hT) -
         \nal(\XPi{\Gz}, \hT)\right) > \sqrt{n}\lambda_n c \notag.
  \end{align}
  As in \mref{eq:normconv}, the left-hand side converges to
  $N(0, \gamma\nu^2)$.

  Furthermore, the sequence $\sqrt{n}\lambda_n$ is bounded and we can apply the
  Bolzano-Weierstrass theorem \cite{bartle}. There must be a subsequence
  $\{\sqrt{n}\lambda_n\}_{n'}$ such that $\lim_{n'\to \infty}
  \sqrt{n'}\lambda_{n'} = \lambda_0 < \infty$. Combining this with
  \mref{eq:toprove} and the asymptotic normality above, we have
  \begin{equation*}
    \nlim \Prob(\SPL{} > \SPL{0}) =
      1 - \Phi\left(\frac{c\lambda_0}{\sqrt{\gamma \nu^2}}\right) > 0,
  \end{equation*}
  where $\Phi$ is standard normal CDF. This proves part (iii).
\end{proof*}

For BIC, $n \lambda_n = \log(n)/2 \to \infty$ and $\sqrt{n}\lambda_n =
\log(n)/(2\sqrt{n}) \to 0$; for AIC, $n \lambda_n = 1$ and
$\sqrt{n}\lambda_n = 1/\sqrt{n} \to 0$. Hence BIC satisfies (i) but not (ii);
and AIC does not satisfy either (i) or (ii), confirming and extending (80) in
\cite{bozdogan}.

Finally, the following corollary justifies the use of NAL in practical BN
structure learning (including CGBNs).

\begin{corollary}[C4.2 in \cite{balov}]
  Assume that \ref{cond:mle}, \ref{cond:secder}, \ref{cond:var} are satisfied.
  For almost all MCAR distributions $\Z$, $[\hG]$ is a consistent estimator of
  $[\Gz]$ if and only if $\lambda_n \to 0$ and $\sqrt{n}\lambda_n\to \infty$.
\label{corConsistency}
\end{corollary}
\begin{proof*}[Corollary \ref{corConsistency}]
  From Corollary \ref{cor:equiv}, we have that if $\X$ is MCAR, $[\Gz]$ is
  identifiable in $\graphs$. Balov \cite{balov} argues that \ref{suff:three}
  holds for almost all distributions of $\Z$ if the set of independence
  relationships implied by $\Gz$ is non-empty. The result then follows from the
  proof of Theorem \ref{theorem:main}, with $\Gz$ replaced by any $\G_* \in
  [\Gz]$.
\end{proof*}

\pagebreak

\section{Computational Experiments}
\label{sec:experiments}

Finally, we assess the structural accuracy and the computational efficiency of
NAL with a simulation study: our aims are to confirm the results presented in
Section \ref{sec:theorems} and to independently validate those originally
presented in Balov \cite{balov}. Firstly, we confirm consistency using exact
learning given a known topological ordering of the nodes and a fixed maximum
number of parents in Section \ref{sec:known-order}. Secondly, we learn DAGs
using tabu search and without conditioning on a topological ordering in
Section \ref{sec:unknown-order}. Besides the accuracy of the learned DAGs, we
explore how different choices of $\lambda_n$ impact the speed of structure
learning. Thirdly, we investigate the claim by Balov \cite{balov} that NAL is an
efficient alternative to Structural EM in Section \ref{sec:em-comparison}.

\begin{table}[b]
\begin{center}
\begin{tabular}{l|l|rr|rr|rr}
             & & & & \multicolumn{2}{l|}{Original} & \multicolumn{2}{l}{Pruned} \\
  Name       & Type & $|\DD|$ & $|\GG|$ & $|A|$ & $|\Theta_0|$ & $|A|$ & $|\Theta_0|$ \\
  \hline
  ALARM      & discrete BN  & 37      & 0       & 46    & 509        & 42    & 239  \\
  ECOLI70    & GBN          & 0       & 46      & 70    & 162        & 62    & 154  \\
  SANGIOVESE & CGBN         & 1       & 14      & 55    & 259        & 35    & 239  \\
  \hline
\end{tabular}
  \caption{Reference BNs from \cite{bnrepo}, with the number of discrete
  ($|\DD|$) and continuous ($|\GG|$) nodes, of arcs ($|A|$) and
  of parameters ($|\Theta_0|$) for the original and the pruned networks.}
\label{tab:ref}
\end{center}
\end{table}

For these purposes, we use three reference BNs from the Bayesian network
repository \cite{bnrepo}: one discrete BN (ALARM), one GBN (ECOLI70) and one
CGBN (SANGIOVESE). Their characteristics are summarised in Table \ref{tab:ref}.
Of the 14 continuous nodes of SANGIOVESE, 4 follow a conditional Gaussian
distribution, while the remaining 10 nodes do not have any discrete parent and
follow a Gaussian distribution.

Our experimental setup largely follows \cite{ijar19}. For computational reasons,
we apply exact learning to pruned versions of the BNs in which each node has at
most two (ALARM, ECOLI) or three (SANGIOVESE) parents. The parameters of pruned
networks are set to the MLEs computed from a data set of size $100 |\Theta_0|$
generated from the original BNs. We generate 20 datasets from each network for
each combination of $n$ (the sample size) and $\beta$ (the probability of a
value to be missing). To facilitate comparisons across networks of different
size and complexity, we choose the values of $n$ as $n = k |\Theta_0|$ to make
them proportional their number of parameters. From each data set, we learn a DAG
for each $\lambda_n$ and we penalise complexity using $h(\G) = |\Theta_\G|$. In
the experiments for known and unknown node ordering we consider all combinations
of the values below, resulting in 144 configurations:
\vspace{-0.5em}
\begin{itemize} \setlength{\itemsep}{-0.25em}
  \item Relative sample size:
    $k \in \{10,\; 50,\; 100,\; 250,\; 500,\; 1000\}$.
  \item Proportion of missing values: $\beta \in \{0,\; 0.05,\; 0.1,\; 0.2\}$.
  \item Penalisation coefficient $\lambda_n$:
    \begin{enumerate}
      \item Custom: $\lambda_n = \frac{1}{N}n^{-\alpha}$ with
        $\alpha \in \{0.10, 0.25, 0.40, 0.60\}$;
      \item BIC: $\lambda_n=\frac{\log n}{2n}$ and AIC:
        $\lambda_n = \frac{1}{n}$.
    \end{enumerate}
\end{itemize}
When comparing NAL with Structural EM we do not consider $k = 1000$ and
$\beta = 0$ and only use $\alpha \in \{0.10, 0.25, 0.40\}$ for computational
reasons.

We compare DAGs using a scaled version of the Structural Hamming Distance (SHD)
by \cite{mmhc}. The SHD between $\G_1$ and $\G_2$ is the minimum number of arc
additions, deletions and reversals required to get from a DAG in $[\G_1]$ to a
DAG in $[\G_2]$. In particular, we use the scaled SHD defined as SHD/$|A|$ in
\cite{ijar19} to allow for meaningful comparisons between networks of different
size. We prefer this measure to other structural distances such as the F1 score
used by Balov \cite{balov} because it takes into account arc directions up to
Markov equivalence. Furthermore, we prefer it to distances in the probability
space such as the Kullback-Leibler divergence because we are mainly interested
in structural accuracy.

We compare NAL and Structural EM on the basis of their running time and of the
number of calls to the scoring function. To make their implementations
comparable, we use the \emph{bnlearn} R package \cite{jss09} to impute missing
data by likelihood-weighted sampling ($500$ particles). For the M-step, we use
BIC as the network score. All the code used in the simulation study is available
from \mbox{\url{https://github.com/Tbodewes/Dissertation}.}

In the following, we discuss results directly in terms of $\alpha$ for
brevity instead of referring to the corresponding $\lambda_n$.

\subsection{Performance for Known Node Ordering}
\label{sec:known-order}

\begin{figure}[t]
  \includegraphics[width=\textwidth]{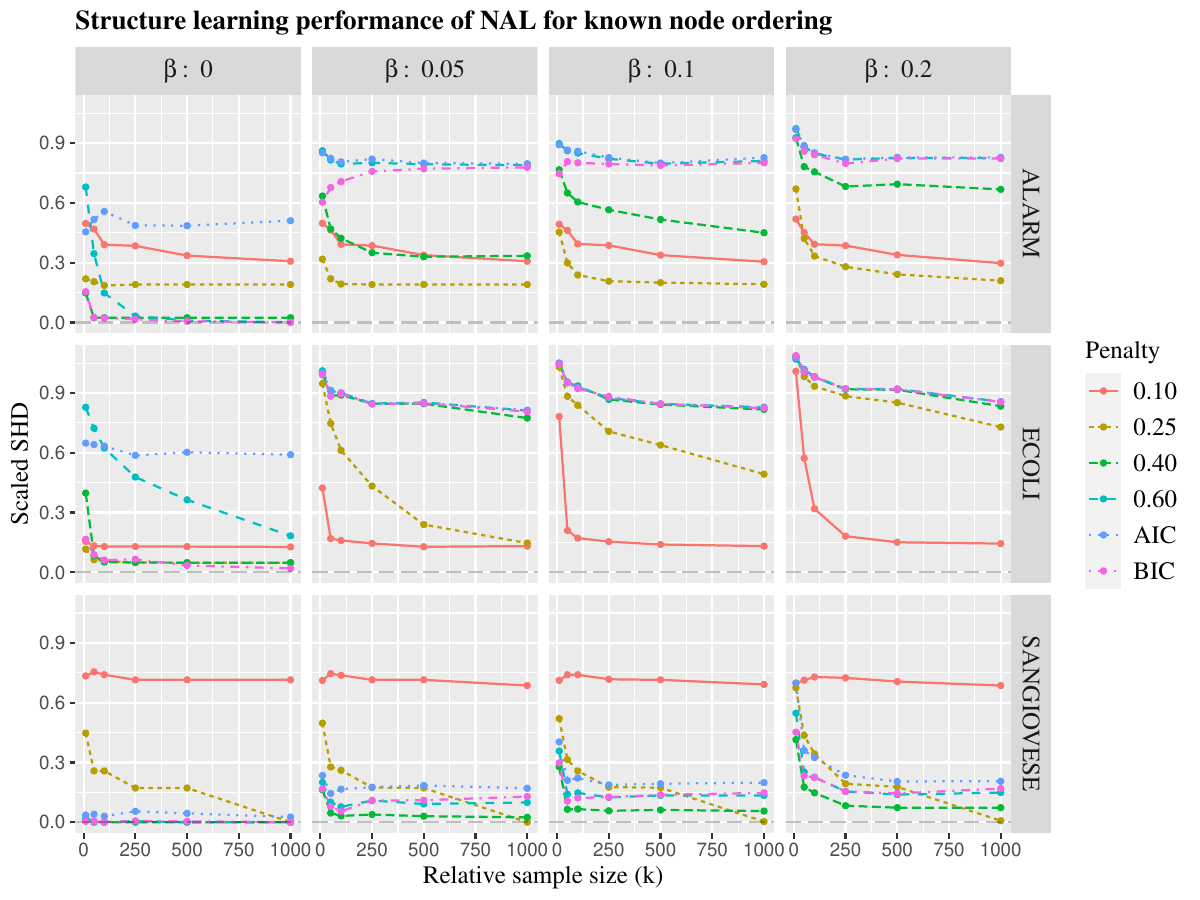}
  \caption{Scaled SHD between DAGs learned using NAL with known node ordering
    and the true DAG. Lower SHD indicates better learning performance. Lines
    correspond to different penalties  $\alpha\in\{0.10, 0.25, 0.40, 0.60\}$,
    BIC and AIC. Error bars are on the order of point size and not shown.}
\label{fig:known-order}
\end{figure}

Figure \ref{fig:known-order} supports our theoretical results and confirms
Balov's \cite{balov}: penalties going to zero slower than $O(n^{-1/2})$ often
lead to increased structural accuracy. In the complete data case, the results
are as expected: BIC is consistent, and penalties other than AIC perform well.
However, $\alpha = 0.10, 0.25$ do not always improve accuracy as the sample size
increases, leading instead to mild underfitting (that is, too few arcs are
included in the DAGs). AIC performs poorly for the opposite reason: it does not
penalise network complexity enough, leading to overfitting (that is, the learned
DAGs contain too many arcs including a number of false positives).

As for incomplete data, we can see that all $\alpha < 0.60$ outperform $\alpha =
0.60$, AIC and BIC for the ALARM network. As the proportion of missing data
increases, the accuracy produced by $\alpha = 0.40$ gradually deteriorates for
all the sample sizes under consideration. We observe similar results for ECOLI:
only $\alpha = 0.10$ does not degrade performance as $\beta$ increases.
Furthermore, $\alpha = 0.25$ does well when $\beta = 0.05$ but it is not
markedly better than the inconsistent scores when $\beta = 0.2$. Closer
inspection suggests that this may be due to overfitting: $\alpha = 0.25, 0.40$
do not penalise complexity enough. Note that Gaussian parents are penalized less
harshly than discrete parents in a CGBN because they introduce fewer parameters
in the BN, which may explain why overfitting occurs for the ECOLI network. Based
on the theoretical results, we expect that penalties below $\alpha = 0.50$ to
give consistent estimates. However, convergence appears to be very slow for
penalties that converge to zero too quickly. For SANGIOVESE, the only penalty
that allows one to learn the true DAG correctly for $k \leqslant 1000$ is $\alpha =
0.25$, with all the $\alpha$s slowly increasing the scaled SHD as $\beta$
increases. $\alpha = 0.10$ seems to suffer from mild underfitting, while $\alpha
= 0.40$, $0.60$, AIC and BIC lead to mild overfitting.

\subsection{Performance for Unknown Node Ordering}
\label{sec:unknown-order}

\begin{figure}[t]
  \includegraphics[width=\textwidth]{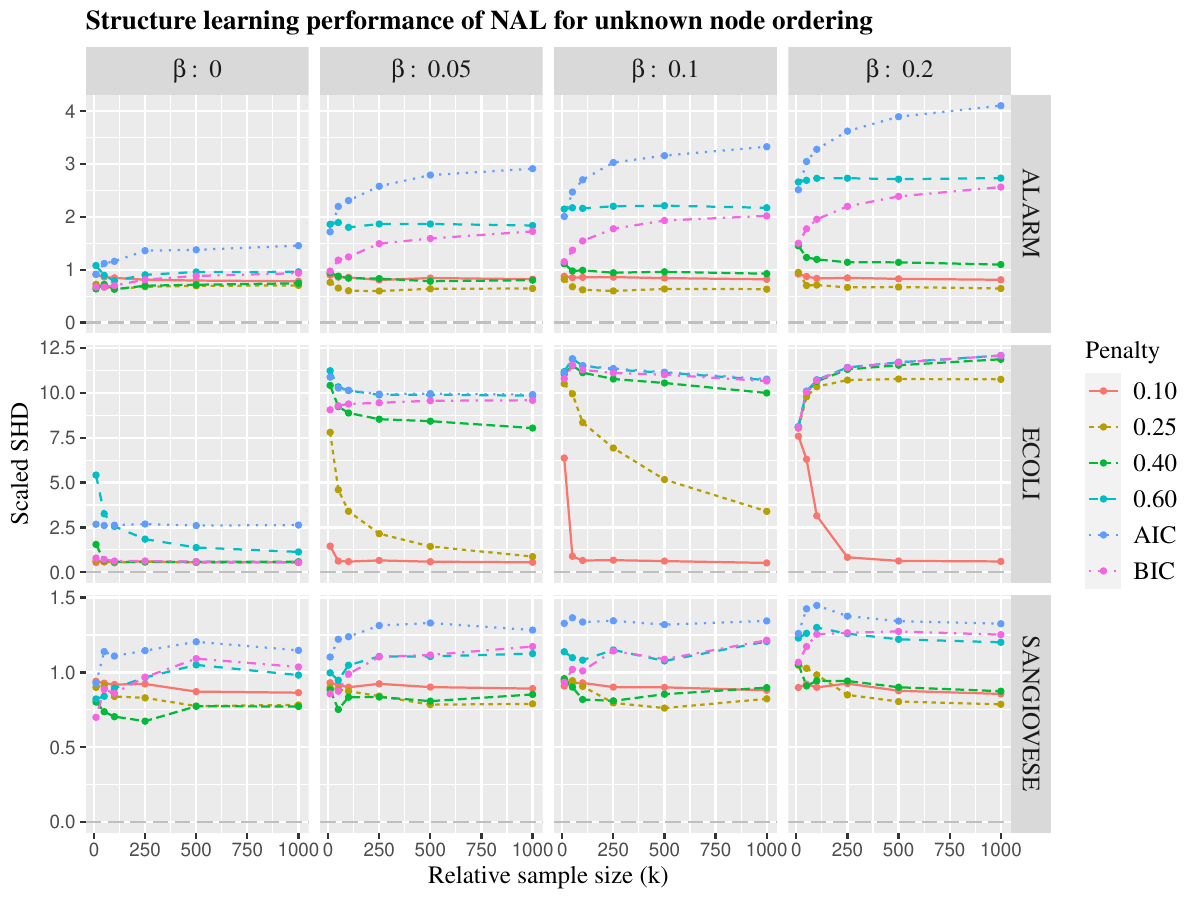}
  \caption{Scaled SHD between DAGs learned using NAL with unknown node ordering
    and true DAG. Note that the scale of the vertical axis differs across rows.
    Error bars are on the order of point size and not shown. For details, see
    the caption of Figure \ref{fig:known-order}.}
  \label{fig:unknown-order}
\end{figure}

Figure \ref{fig:unknown-order} suggests that for unknown node ordering,
$\alpha = 0.10, 0.25$ tend to give the best results. In the complete data case,
BIC no longer outperforms $\alpha = 0.10, 0.25$ but the overall behaviour of
other penalties is similar to that shown in Figure \ref{fig:known-order}.
All penalties except AIC give a similar scaled SHD for each BN. Penalties that
allowed one to learn the true DAG when the topological ordering of the nodes was
known now have a scaled SHD strictly greater than zero due to the heuristic
nature of tabu search: SHD/$|A| \approx 1$ for all of them. This is in line with
the findings in \cite{ijar19}.

As for incomplete data, we see a clean separation between consistent and
inconsistent scores for the ALARM network. As $\beta$ increases, the performance
of AIC, BIC and $\alpha = 0.60$ gradually degrades, while the performance of
$\alpha = 0.10, 0.25$ remains broadly the same. For ECOLI, only $\alpha = 0.10$
results in scaled SHD values that are promising for practical applications, with
the other penalties giving scaled SHDs around or larger 10. (That is, each
DAG we learned includes 10 times as many arcs as the corresponding true DAG.)
The scaled SHD of $\alpha = 0.25$ is close to that of $\alpha = 0.10$ for $\beta
= 0.05$, but increases rapidly as $\beta$ grows. This is due to the overfitting
also observed for known node ordering, but in this case it is more pronounced
because we did not prune parent sets. For SANGIOVESE, SHD/$|A| \approx 1$ for
all penalties. Nevertheless, performance is better when $\lambda_n$ goes to zero
at a rate slower than $n^{-1/2}$.

\begin{figure}[t]
  \includegraphics[width=\textwidth]{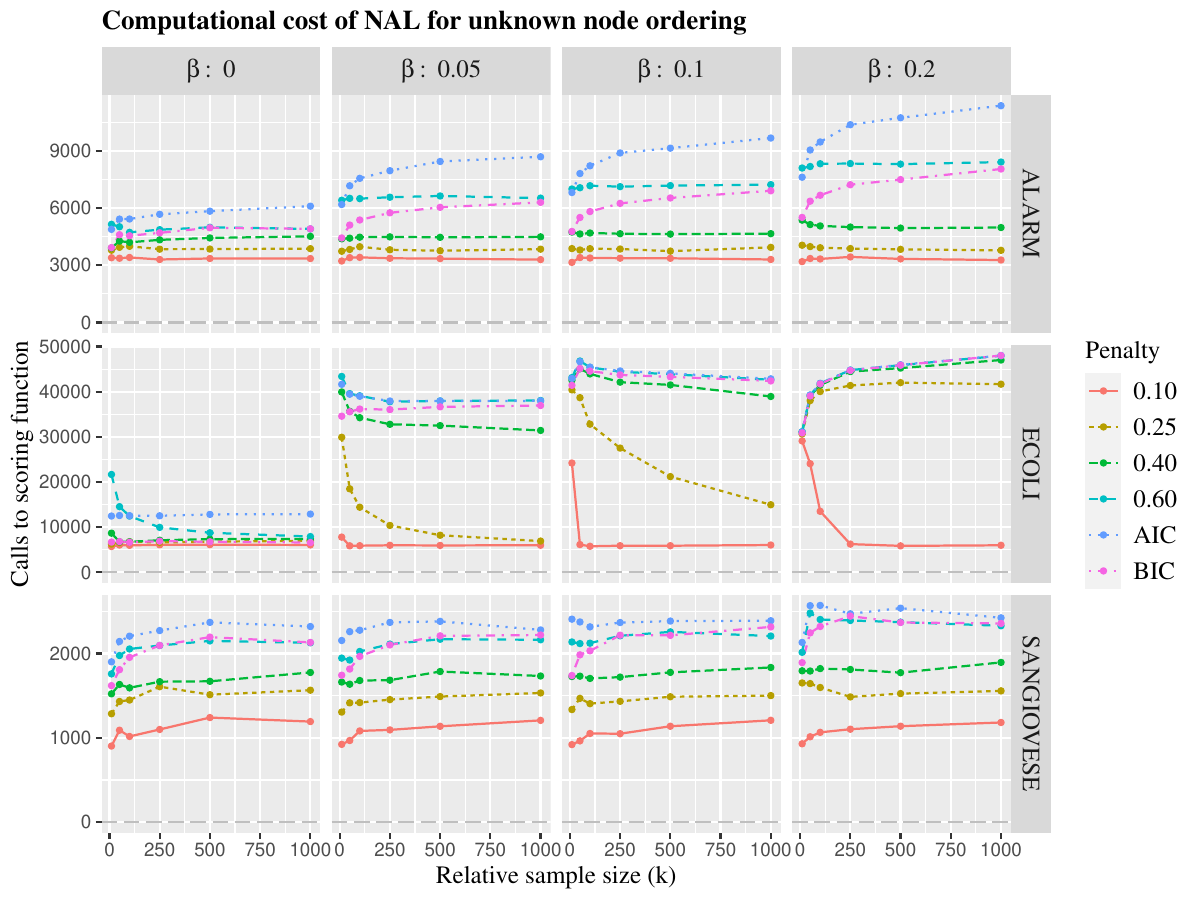}
  \caption{Number of calls to the scoring function until tabu search converged.
    Note that the scale of the vertical axis differs across rows. Error bars are
    on the order of graph point size and not shown. For details, see the caption
    of Figure \ref{fig:known-order}.}
  \label{fig:queries-unknown}
\end{figure}

Figure \ref{fig:queries-unknown} shows that stronger penalties generally lead
to tabu search converging faster. For ALARM, $\alpha = 0.10$ leads to a
50--70\% reduction in the number of calls to the scoring function relative to
$\alpha = 0.40$, with even greater reductions for the other BNs. The
computational cost follows a pattern that is similar to the SHD pattern in
Figure \ref{fig:unknown-order}: inconsistent scores have over twice the
computation cost of the consistent scores because weaker penalties lead to more
arcs being included in the learned DAGs.

\subsection{Comparison of Node-Average Likelihood and Structural EM}
\label{sec:em-comparison}

\begin{figure}[t]
  \includegraphics[width=\textwidth]{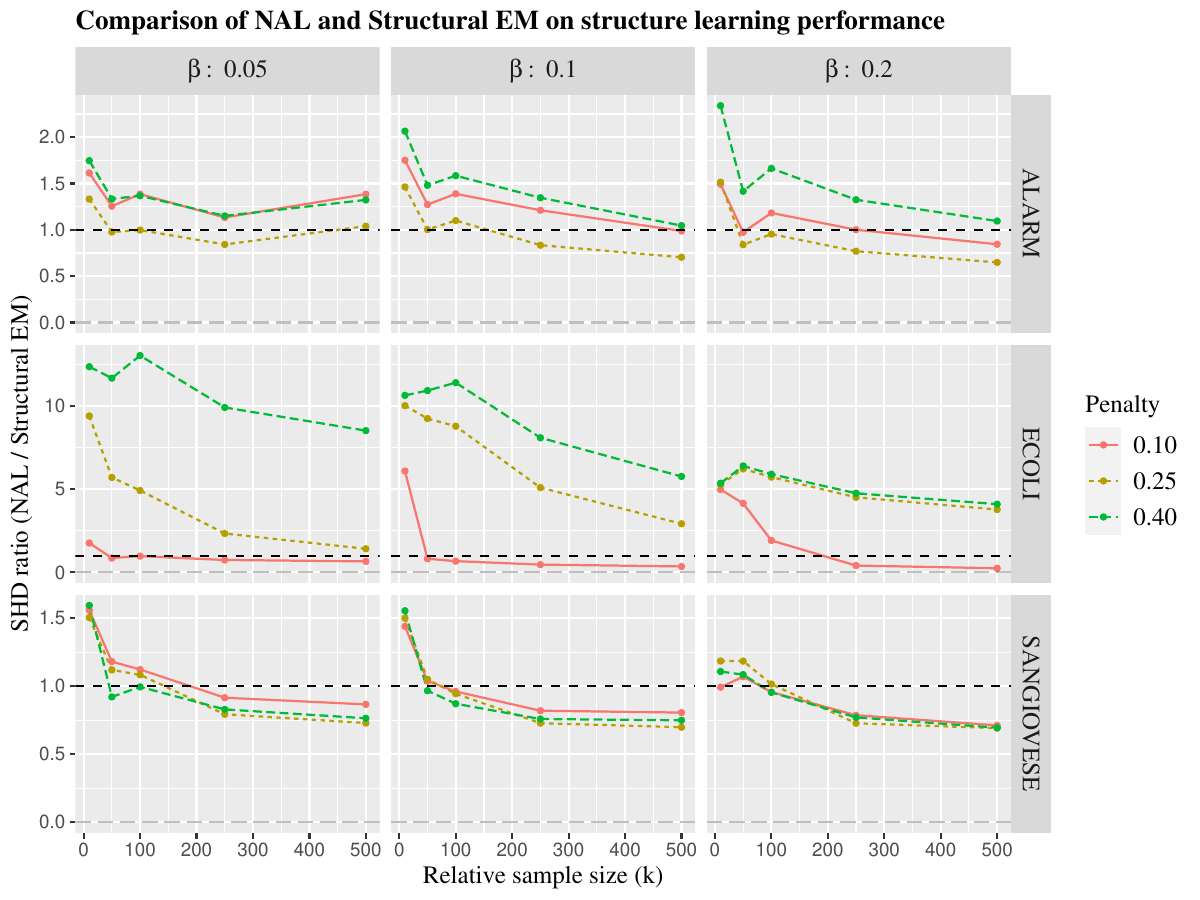}
  \caption{Ratio between the SHD for the DAG learned using NAL and that
  for the DAG learned using Structural EM, both with respect to the true DAG.
  The black dashed line marks equal performance of NAL and Structural EM. A
  lower ratio indicates better relative performance of NAL. Lines correspond to
  different penalties $\alpha \in \{0.10, 0.25, 0.40\}$. Error bars are on the
  order of point size and not shown.}
  \label{fig:em-shd}
\end{figure}

Figure \ref{fig:em-shd} suggests that NAL and Structural EM have similar
structural accuracy, with NAL improving relative to Structural EM for larger
values of $k$. For ALARM, NAL outperforms Structural EM for larger $k$ when
$\beta \geqslant 0.1$ and $\alpha = 0.25$. For ECOLI, only NAL with $\alpha =
0.10$ is comparable to or better than Structural EM in terms of SHD. This is in
line with the performance gap observed between $\alpha = 0.10$ and all other
penalties in Figure \ref{fig:unknown-order}, and suggests that Structural EM
does not suffer as strongly as NAL from the overfitting observed in Figure
\ref{fig:unknown-order}. For SANGIOVESE, the two methods perform similarly
again, with NAL outperforming Structural EM for larger sample sizes.

\begin{figure}
\includegraphics[width=\linewidth]{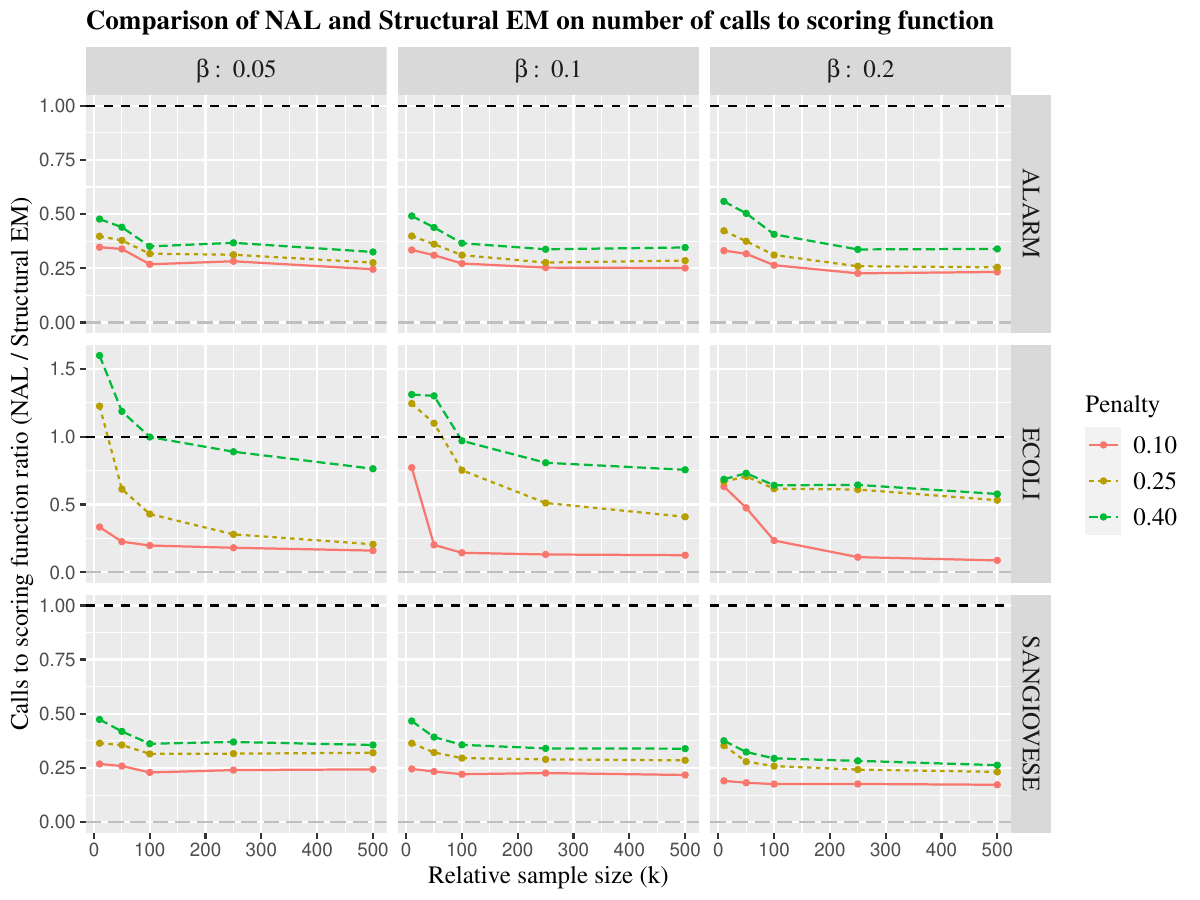}
\includegraphics[width=\linewidth]{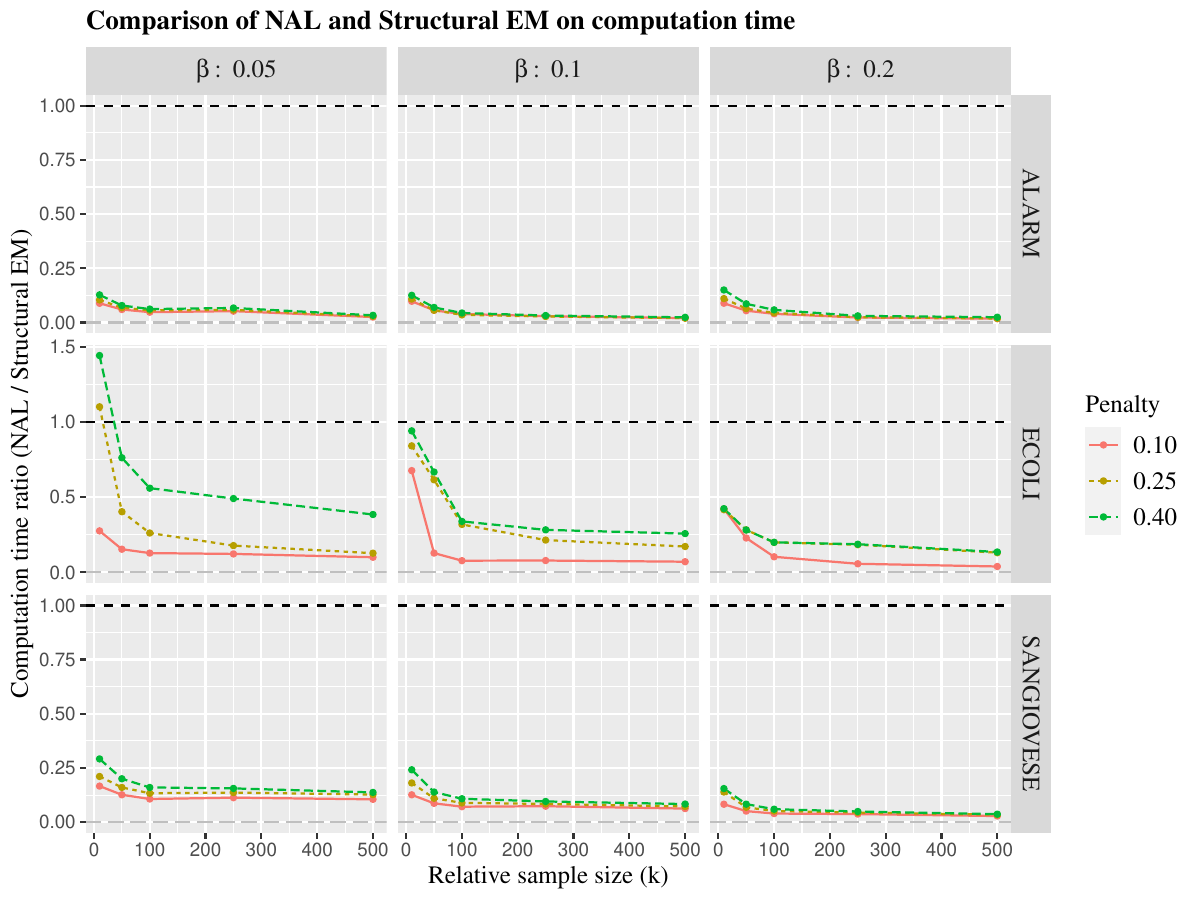}
\caption{Comparison of NAL and Structural EM in terms of number of calls to the
  scoring function (top panel) and single-core running time (bottom panel).
  Formatting is the same as in Figure \ref{fig:em-shd}.}
\label{fig:em-timing}
\end{figure}

Figure \ref{fig:em-timing} (top panel) shows that for appropriate penalisation,
NAL requires 50--75\% fewer calls to the scoring function than Structural EM. For
ALARM and SANGIOVESE this holds for all considered penalties and sample sizes.
For ECOLI, this is only true for larger sample sizes or for $\alpha = 0.10$.
This agrees with the evidence from Figure \ref{fig:queries-unknown} that NAL
with $\alpha = 0.25,0.40$ converges slowly for ECOLI. Figure \ref{fig:em-timing}
(bottom panel) shows that NAL is four to ten times faster than Structural EM.
The difference is larger than the difference in the number of queries due to the
computational cost of the E-step in Structural EM.

\section{Conclusions}

Common approaches to BN structure learning from incomplete data, such as the
Structural EM, embed EM within score-based algorithms thus incurring a
significant computational cost. NAL is a competitive but faster alternative; in
this paper we proved its consistency and model identifiability for CGBNs,
showing NAL's wide applicability beyond its original formulation. We also
established general sufficient conditions for consistency that can be readily
checked for other classes of BNs. To confirm these theoretical results and
Balov's \cite{balov}, as well as Balov's assertion that NAL is faster than
Structural EM, we performed a simulation study including a discrete BN, a GBN
and a CGBN. We observed that indeed:
\begin{itemize}
  \item BIC is not consistent with incomplete data, and AIC is not consistent
    even with complete data.
  \item NAL with appropriate penalties has similar structural accuracy to
    Structural EM.
  \item NAL with appropriate penalties can require 50-75\% fewer calls to the
    scoring function than Structural EM.
  \item NAL with appropriate penalties can converge four to ten times faster
    than Structural EM.
\end{itemize}
These observations suggest that NAL is competitive with Structural EM in terms
of structural accuracy, while being computationally more efficient, in the
broader class of CGBN rather than just for discrete BNs.

\section*{Acknowledgements}

This work is an extended version of the paper ``Identifiability and Consistency
of Bayesian Network Structure Learning from Incomplete Data'' published in the
Proceedings of the 10th International Conference on Probabilistic Graphical
Models \cite{pgm20-tjebbe}.

\appendix

\section{Proofs of Auxiliary Results}
\label{sec:proofs}

\begin{proof*}[Lemma \ref{lemma:nondec}]
  Let $\ZA = \Z_{X_i, \A}$ and $\ZAB = \Z_{X_i, \A, \B}$. By the law of total
  probability and $\Z_\B \indep X_i, \B \given \Z_{(\A)}$:
  \begin{equation*}
    \Prob(X_i \given \A, \ZA) =
      \E\left[\Prob(X_i \given \A, \B, \ZAB) \given \A, \ZAB\right].
  \end{equation*}
  Since $\Z_\B \indep \A \given \ZA$, by Jensen's inequality
  \begin{align*}
    \bar\logl(X_i \given \A)
    &= \E\left[\int \Prob(X_i \given \A, \ZA)
         \log\Prob(X_i \given \A,  \ZA)\dd X_i \at \ZA\right] \\
    &\leqslant \E\left[ \int \Prob(X_i \given \A, \B, \ZAB)
      \log\Prob(X_i \given \A, \B, \ZAB) \dd X_i \at \ZAB\right] \\
    &= \bar\logl(X_i \given \A, \B),
  \end{align*}
  with equality if and only if $X_i \indep \B \given (\A, \ZA)$.
\end{proof*}
\begin{proof*}[Proposition \ref{prop:id}]
  Under MCAR, $\Prob_{\Gz}(\X \given \Z) = \Prob_{\Gz}(\X)$ and
  \begin{equation*}
    \nal(\G, \Theta) \approx \E\left[ \nal(\G, \Theta) \right] =
    \E\left[ \logl(\Gz, \Tz) \right] - \KL(\Prob_{\Gz}(\X) \klsep \Prob_\G(\X)).
  \end{equation*}
  The Kullback-Leibler divergence $\KL(\Prob_{\Gz}(\X) \klsep \Prob_\G(\X)) = 0$
  if and only if $\Prob_{\Gz}(\X) = \Prob_\G(\X)$ almost everywhere, hence
  $\max_{\G \inG} \logl(\G, \Theta) = \logl(\Gz, \Tz)$. The maximum is attained
  only when $\Prob_{\Gz}(\X) = \Prob_\G(\X)$. The maximum can always be attained
  because $\Gz \inG$. This proves \ref{part:one} and \ref{part:two}.
  For \ref{part:three}, $\PXi{\Gz} \subseteq \PXi{\G}$ for all $X_i$ and
  due to Lemma \ref{lemma:nondec} $\nal(\G, \Theta) \geqslant \nal(\Gz, \Tz)$,
  while \ref{part:one} implies $\nal(\G, \Theta) \leqslant \nal(\Gz, \Tz)$.
  $\nal(\G, \Theta) = \nal(\Gz, \Tz)$ follows.
  Not that $\nal(\G, \Theta)$ converges to its expectation as $n \to \infty$,
  so the approximation above becomes an equality in the limit.
\end{proof*}
\begin{proof*}[Corollary \ref{cor:equiv}]
  Let $\G_*$ satisfy $\nal(\G_*, \Theta_*) = \max_{\G \inG} \nal(\G, \Theta)$.
  By \ref{part:one} of Proposition \ref{prop:id}, $\nal(\G_*, \Theta_*) =
  \nal(\Gz, \Tz)$. Then, by \ref{part:two} we have $\Prob_{\G}(\X) =
  \Prob_{\Gz}(\X)$. The minimum of all such $\G_*$ is a valid DAG. The result
  follows from \cite{chickering}.
\end{proof*}
\begin{proof*}[Proposition \ref{prop:regCGBN}]
  \textit{Condition \ref{cond:mle}}: Let $\Ij$ be the subset of observations for
  which $\DXi = j$. By assumption, $P(\DXi = j) > 0$, thus this set is non-empty
  for $n \to \infty$. If $X_i \in \DD$, $\hT$ satisfies $\Jac \nal(\XPi{}, \hT)
  = 0$ and is given by $\hpiijk = \dim{\Ij}^{-1}\sum_{\Ij} \bbbone[X_i = k]$.
  The Lindeberg-L\'evy CLT gives $\hpiijk = \piijk + O(n^{-1/2})$ as desired.

  Now consider $X_i \in \GG$. From \mref{eq:cgnode}, if we let $\W =
  [\mathbf{1} \; \GXi]$ conditional on $\Ij$ and we absorb $\mXd$ into $\bXd$,
  the MLEs for $\bXd$ and $\sXd$ are the classic
  \begin{align*}
    \hbXd &= (\W^T \W)^{-1}\W^T X_i= \beta_k + (\W^T \W)^{-1} \W^T \varepsilon, \\
    \hsXd &= \frac{\widehat{e}^T \widehat{e}}{n} =
      \frac{\varepsilon^T \varepsilon}{n} - \frac{1}{n}
      \left(\frac{\W^T \varepsilon}{\sqrt{n}}\right)
      \left(\frac{\W^T \W}{n}\right)^{-1}
      \left(\frac{\W^T \varepsilon}{\sqrt{n}}\right).
  \end{align*}
  with $\widehat{e} = X_i - \W^T\hbXd$. $(\W^T \W/n)^{-1} \to Q^{-1}$ and $Q$ is
  invertible under our assumptions. Hence $\hbXd$ and $\hsXd$ satisfy
  $\Jac \nal(\XPi{}, \hT) = 0$, and as $\W$ and $\varepsilon$ are uncorrelated
  \begin{equation*}
    \sqrt{n}(\hbXd - \bXd) =
      \left(\W^T \W/n \right)^{-1}\left(\W^T\varepsilon/\sqrt{n}\right) = O(1)
  \end{equation*}
  and
  \begin{equation*}
    \sqrt{n}\left(\hsXd - \sXd\right) = O(1)
  \end{equation*}
  as shown, for instance, in \cite{heij}. \ref{cond:mle} follows.

  \textit{Condition \ref{cond:secder}}: For $X_i \in \DD$, the diagonal elements
  of the Hessian are $-1/\pi_{jk}^2$, while the off-diagonal elements are zero.
  Hence the Hessian is finite for $\pi_{jk}>0$.

  For $X_i \in \GG$, the Hessian of $\bXd, \sXd$ for each $\DXi = j$ is
  \begin{equation*}
    -\left(\sXd\right)^{-2}\begin{bmatrix}
      \sXd \W\W^T & e \W^T\\
      e^T \W & (2e^2 - \sXd)/(2\sXd)
    \end{bmatrix},
  \end{equation*}
  where $e = X_i - \W^T \bXd$. All elements of this matrix are constants,
  Gaussian random variables or squares and cross-products thereof. As $\sXd >0$
  by assumption, the expectation of all elements of the Hessian is finite and
  \ref{cond:secder} is satisfied.

  \textit{Condition \ref{cond:var}}:
  By Jensen's inequality, it suffices to show that \newline
  $\E[\nal(\XPi{}, \T)^2] < \infty$. If $X_i \in \DD$,
  \begin{equation*}
    \E\left[\nal(\XPi{}, \T)^2\right] =
      \sum\nolimits_{j, k} \piijk (\log \piijk)^2 < \infty
  \end{equation*}
  as we are summing a finite number of finite terms under the assumption
  that all $\piijk > 0$. For $X_i \in \GG$,
  \begin{equation*}
    \E\left[\nal(\XPi{}, \T)^2\right] =
      \sum\nolimits_j \Prob(\DXi = j) \E\left[(\log\Prob(X_i \given \GXi,
        \DXi = j, \bXd, \sXd))^2\right]
  \end{equation*}
  where the expectations on the right-hand side, disregarding constants, take
  the form
  \begin{equation*}
    \E\left[\left((X_i - \GXi\bXd)^2/(2\sXd)\right)^2\right] < \infty
  \end{equation*}
  since Gaussian moments are finite. Thus \ref{cond:var} is satisfied.
\end{proof*}
\begin{proof*}[Lemma \ref{lemma:toscore}]
  The difference in scores, scaled by $n^{\alpha}$, is
  \begin{multline*}
    n^{\alpha}(\SPL{1} - \SPL{2}) = \\
    n^{\alpha}(\nal(\G_1, \Theta_1 \given \D) - \nal(\G_2, \Theta_2 \given \D)) +
    n^{\alpha}\lambda_n (h(\G_2) - h(\G_1)).
  \end{multline*}
  The first term is $O(1)$, while the second diverges as $n \to \infty$. The
  result follows.
\end{proof*}
\begin{proof*}[Lemma \ref{lemma:mle}]
  If $n \to \infty$, under \ref{cond:mle} and \ref{cond:secder}
  \begin{multline*}
    \nal(\XPi{}, \T) = \nal(\XPi{}, \hT) + \\
      \frac{1}{2}(\T - \hT)^T \Hess\left(\nal(\XPi{}, \T)\right)(\T - \hT).
  \end{multline*}
  All elements of the Hessian are bounded by \ref{cond:secder}, thus
  \begin{equation*}
    \nal(\XPi{}, \hT) - \nal(\XPi{}, \T) = O(\lVert \T - \hT \rVert^2),
  \end{equation*}
  which is $O(n^{-1})$. The result follows.
\end{proof*}
\begin{proof*}[Lemma \ref{lemma:cpd}]
  Under MCAR, we can drop $\Zi$. If $\G = \Gz$, (i) and (ii) are
  trivial, while (iii) holds by definition. If $\G \supset \Gz$, there exists an
  $X_i$ such that $\PXi{\Gz} \subset \PXi{\G}$ and for which
  $X_i \indep \PXi{\G} \setminus \PXi{\Gz} \given \PXi{\Gz}$,
  then (i) follows directly. For (ii), note that
  $\Prob(\XPi{\G}) = \Prob(\XPi{\G^0})$.
  By definition, estimating $\T$ minimises the $\KL$ divergence between $\G,
  \Theta$ and $\Gz, \Tz$ for the given $\D$, which happens when the
  corresponding distributions are equal almost surely leading to (ii). We obtain
  (iii) from (ii) as $\Prob(\XPi{\G}, \T, \D) = \Prob(\XPi{\G^0}, \T, \D)$
  \begin{equation*}
    \Prob(\XPi{\G^0}, \T, \D) =
      \Prob(\XPi{\G^0}, \T) = \Prob(\XPi{\G^0}, \T, \Zi).
  \end{equation*}
\end{proof*}
\begin{proof*}[Lemma \ref{lemma:lconv}]
  Under MCAR, we can drop the $\Zi$. Given \mref{eq:loglik} and the regularity
  conditions \ref{cond:mle}, \ref{cond:secder}, Lemma \ref{lemma:mle} gives
  \begin{equation}
    \nal(\XPi{}, \hT) = \nal(\XPi{}, \T) + O(n^{-1}).
  \label{eq:substitution}
  \end{equation}
  As $\nal(\XPi{}, \T)$ is a sum of i.i.d. random variables whose mean exists
  and is finite, we can apply the Weak Law of Large Numbers to the right hand
  side of \mref{eq:substitution}. Then
  \begin{multline*}
    \nal(\XPi{}, \hT) \approx \nal(\XPi{}, \T) - \\
       \E\left[\KL\left(\Prob(\XPi{}, \T) \klsep \Prob(\XPi, \hT))\right)\right].
  \end{multline*}
  $\KL(\cdot) \geqslant 0$ implies that
  $\nlim \Prob(\nal(\XPi{}, \hT) \leqslant \nal(\XPi{}, \T)) = 1$
  and \mref{eq:negativity} follows.

  To prove the second assertion, assume that $\Gz \subseteq \G$. By Lemma
  \ref{lemma:cpd}(iii), we have $\Prob(\XPi{}, \T, \Zi) = \Prob(\XPi{}, \T, \D)$
  almost surely. Hence the KL divergence above is zero and
  $\nal(\XPi{}, \hT) \to_p \nal(\XPi{}, \T)$
  with a rate of convergence of $O(n^{-1/2})$ by the Lindeberg-L\'evy CLT,
  dominating the $O(n^{-1})$ in \mref{eq:substitution}.
\end{proof*}


\begin{thebibliography}{10}
\expandafter\ifx\csname url\endcsname\relax
  \def\url#1{\texttt{#1}}\fi
\expandafter\ifx\csname urlprefix\endcsname\relax\def\urlprefix{URL }\fi
\expandafter\ifx\csname href\endcsname\relax
  \def\href#1#2{#2} \def\path#1{#1}\fi

\bibitem{balov}
N.~Balov, {Consistent Model Selection of Discrete Bayesian Networks from
  Incomplete Data}, Electronic Journal of Statistics 7 (2013) 1047--1077.

\bibitem{pearl}
J.~Pearl, {Probabilistic Reasoning in Intelligent Systems: Networks of
  Plausible Inference}, Morgan Kaufmann Publishers Inc., 1988.

\bibitem{lauritzen}
S.~L. Lauritzen, {Graphical Models}, Oxford University Press, 1996.

\bibitem{chickering}
D.~M. Chickering, {A Transformational Characterization of Equivalent Bayesian
  Network Structures}, in: {UAI}, 1995, pp. 87--98.

\bibitem{heckerman}
D.~Heckerman, D.~Geiger, D.~M. Chickering, {Learning Bayesian Networks: The
  Combination of Knowledge and Statistical Data}, Machine Learning 20~(3)
  (1995) 197--243.

\bibitem{heckerman3}
D.~Geiger, D.~Heckerman, {Learning Gaussian Networks}, in: {UAI}, 1994, pp.
  235--243.

\bibitem{weatherburn}
C.~E. Weatherburn, {A First Course in Mathematical Statistics}, Cambridge
  University Press, 1961.

\bibitem{unification}
D.~Heckerman, D.~Geiger, {Learning Bayesian Networks: a Unification for
  Discrete and Gaussian Domains}, in: UAI, 1995, pp. 274--284.

\bibitem{copula}
G.~Elidan, {Copula Bayesian Networks}, in: J.~D. Lafferty, C.~K.~I. Williams,
  J.~Shawe-Taylor, R.~S. Zemel, A.~Culotta (Eds.), Advances in Neural
  Information Processing Systems 23, 2010, pp. 559--567.

\bibitem{truncexp}
S.~Moral, R.~Rumi, A.~Salmer{\'o}n, {Mixtures of Truncated Exponentials in
  Hybrid Bayesian Networks}, in: {Symbolic and Quantitative Approaches to
  Reasoning with Uncertainty (ECSQARU)}, Vol. 2143 of Lecture Notes in Computer
  Science, Springer, 2001, pp. 156--167.

\bibitem{nphard}
D.~M. Chickering, D.~Heckerman, {Learning Bayesian Networks is NP-hard}, Tech.
  Rep. MSR-TR-94-17, Microsoft Corporation (1994).

\bibitem{ic}
T.~S. Verma, J.~Pearl, {Equivalence and Synthesis of Causal Models},
  Uncertainty in Artificial Intelligence 6 (1991) 255--268.

\bibitem{colombo}
D.~Colombo, M.~H. Maathuis, {Order-Independent Constraint-Based Causal
  Structure Learning}, Journal of Machine Learning Research 15 (2014)
  3921--3962.

\bibitem{spirtes}
P.~Spirtes, C.~Glymour, R.~Scheines, {Causation, Prediction, and Search}, MIT
  Press, 2000.

\bibitem{bouckaert}
R.~R. Bouckaert, {Bayesian Belief Networks: from Construction to Inference},
  Ph.D. thesis, Utrecht University, The Netherlands (1995).

\bibitem{larranaga}
P.~Larra{\~n}aga, B.~Sierra, M.~J. Gallego, M.~J. Michelena, J.~M. Picaza,
  {Learning Bayesian Networks by Genetic Algorithms: A Case Study in the
  Prediction of Survival in Malignant Skin Melanoma}, in: {Proceedings of the
  6th Conference on Artificial Intelligence in Medicine in Europe (AIME'97)},
  Springer, 1997, pp. 261--272.

\bibitem{order1}
P.~Larra{\~n}aga, C.~M.~H. Kuijpers, R.~H. Murga, Y.~Yurramendi, {Learning
  Bayesian Network Structures by Searching for the Best Ordering with Genetic
  Algorithms}, {IEEE Transactions on Systems, Man, and Cybernetics - Part A:
  Systems and Humans} 26~(4) (1996) 487--493.

\bibitem{k2}
G.~F. Cooper, E.~Herskovits, {A Bayesian Method for Constructing Bayesian
  Belief Networks from Databases}, in: {Proceedings of the 7th Conference on
  Uncertainty in Artificial Intelligence}, 1991, pp. 86--94.

\bibitem{cutting}
J.~Cussens, {Bayesian Network Learning with Cutting Planes}, in: {Proceedings
  of the 27th Conference on Uncertainty in Artificial Intelligence}, 2012, pp.
  153--160.

\bibitem{suzuki17}
J.~Suzuki, {An Efficient Bayesian Network Structure Learning Strategy}, New
  Generation Computing 35~(1) (2017) 105--124.

\bibitem{mmhc}
I.~Tsamardinos, L.~E. Brown, C.~F. Aliferis, {The Max-Min Hill-Climbing
  Bayesian Network Structure Learning Algorithm}, Machine Learning 65~(1)
  (2006) 31--78.

\bibitem{ijar19}
M.~Scutari, C.~E. Graafland, J.~M. Guti{\'e}rrez, {Who Learns Better Bayesian
  Network Structures: Accuracy and Speed of Structure Learning Algorithms},
  International Journal of Approximate Reasoning 115 (2019) 235--253.

\bibitem{em}
A.~P. Dempster, N.~M. Laird, D.~B. Rubin, {Maximum Likelihood from Incomplete
  Data via the EM Algorithm}, Journal of the Royal Statistical Society, Series
  B (1977) 1--38.

\bibitem{friedman97}
N.~Friedman, {Learning Belief Networks in the Presence of Missing Values and
  Hidden Variables}, in: ICML, 1997, pp. 125--133.

\bibitem{friedman98}
N.~Friedman, {The Bayesian Structural EM Algorithm}, in: UAI, 1998, pp.
  129--138.

\bibitem{neerlandica19}
M.~Scutari, Bayesian network models for incomplete and dynamic data, Statistica
  NeerlandicaIn print. (2020).

\bibitem{glover}
F.~Glover, M.~Laguna, Tabu search, Springer, 1998.

\bibitem{harary}
F.~Harary, E.~M. Palmer, {Graphical Enumeration}, Academic Press, 1973.

\bibitem{pgm16}
M.~Scutari, {An Empirical-Bayes Score for Discrete Bayesian Networks}, Journal
  of Machine Learning Research (Proceedings Track, PGM 2016) 52 (2016)
  438--448.

\bibitem{gprior}
R.~Eggeling, J.~Viinikka, A.~Vuoksenmaa, M.~Koivisto, {On Structure Priors for
  Learning Bayesian Networks}, {Proceedings of Machine Learning Research
  (AISTATS 2019)} 89 (2019) 1687--1695.

\bibitem{csprior}
R.~Castelo, A.~Siebes, {Priors on Network Structures. Biasing the Search for
  Bayesian Networks}, International Journal of Approximate Reasoning 24~(1)
  (2000) 39--57.

\bibitem{mukherjee}
S.~Mukherjee, T.~P. Speed, {Network Inference Using Informative Priors},
  Proceedings of the National Academy of Sciences 105~(38) (2008) 14313--14318.

\bibitem{ambn17}
M.~Scutari, {Dirichlet Bayesian Network Scores and the Maximum Entropy
  Principle}, Proceedings of Machine Learning Research (AMBN 2017) 73 (2017)
  9--20.

\bibitem{schwarz}
G.~Schwarz, {Estimating the Dimension of a Model}, The Annals of Statistics
  6~(2) (1978) 461--464.

\bibitem{l0}
S.~{van de Geer}, P.~B{\"u}hlmann, {$\ell_0$-Penalized Maximum Likelihood for
  Sparse Directed Acyclic Graphs}, The Annals of Statistics 41~(2) (2013)
  536--567.

\bibitem{rissanen}
J.~Rissanen, {Information and Complexity in Statistical Models}, Springer,
  2007.

\bibitem{akaike}
H.~Akaike, {A New Look at the Statistical Model Identification}, IEEE
  Transactions on Automatic Control 19~(6) (1974) 716 -- 723.

\bibitem{bozdogan}
H.~Bozdogan, {Model Selection and Akaike's Information Criterion (AIC): The
  General Theory and its Analytical Extensions}, Psychometrika 52~(3) (1987)
  345--370.

\bibitem{rubin2}
D.~B. Rubin, {Inference and Missing Data}, Biometrika 63 (1976) 581--592.

\bibitem{rubin1}
R.~J.~A. Little, D.~B. Rubin, {Statistical Analysis with Missing Data}, 1st
  Edition, Wiley, 1987.

\bibitem{grem}
S.~L. Lauritzen, {The EM algorithm for Graphical Association Models with
  Missing Data}, Computational Statistics \& Data Analysis 19~(2) (1995)
  191--201.

\bibitem{shafer}
G.~Shafer, P.~P. Shenoy, {Probability propagation}, Annals of Mathematics and
  Artificial Intelligence 2~(1-4) (1990) 327--351.

\bibitem{variational-em}
M.~Beal, Z.~Ghahramani, {The Variational Bayesian EM Algorithm for Incomplete
  Data: with Application to Scoring Graphical Model Structures}, Bayesian
  Statistics 7 (2003) 453--464.

\bibitem{billingsley}
P.~Billingsley, {Probability and Measure}, Wiley, 1995.

\bibitem{glm}
P.~McCullagh, J.~A. Nelder, {Generalized Linear Models}, 2nd Edition, CRC
  press, 1989.

\bibitem{fahrmeir}
L.~Fahrmeir, H.~Kaufmann, {Consistency and Asymptotic Normality of the Maximum
  Likelihood Estimator in Generalized Linear Models}, The Annals of Statistics
  (1985) 342--368.

\bibitem{white}
H.~White, {Maximum Likelihood Estimation of Misspecified Models}, Econometrica
  50~(1) (1982) 1--25.

\bibitem{bartle}
R.~G. Bartle, D.~R. Sherbert, {Introduction to Real Analysis}, Wiley, 2000.

\bibitem{bnrepo}
M.~Scutari, {Bayesian Network Repository}, http://www.bnlearn.com/bnrepository
  (2012).

\bibitem{jss09}
M.~Scutari, {Learning Bayesian Networks with the bnlearn R Package}, Journal of
  Statistical Software 35~(3) (2010) 1--22.

\bibitem{pgm20-tjebbe}
T.~Bodewes, M.~Scutari, {Identifiability and Consistency of Bayesian Network
  Structure Learning from Incomplete Data}, Proceedings of Machine Learning
  Research (PGM 2020) 138 (2020) 29--40.

\bibitem{heij}
C.~Heij, , P.~{de Boer}, P.~H. Franses, T.~Kloek, H.~K. {van Dijk},
  {Econometric Methods with Applications in Business and Economics}, OUP
  Oxford, 2004.

\end{thebibliography}
\end{document}